\documentclass{amsart}

\usepackage{amsmath,amssymb,graphics}
\usepackage[headings]{fullpage}

\makeatletter
\def\ps@headings{\ps@empty
  \def\@evenhead{%
    \setTrue{runhead}%
    \normalfont\scriptsize
    \rlap{\thepage}\hfil
    \def\thanks{\protect\thanks@warning}%
    \leftmark{}{}\hfil}%
  \def\@oddhead{%
    \setTrue{runhead}%
    \normalfont\scriptsize \hfil
    \def\thanks{\protect\thanks@warning}%
    \rightmark{}{}\hfil \llap{\thepage}}%
  \let\@mkboth\markboth
}
\makeatother
  \def\cal{\mathcal}

\usepackage[all]{xy}

\newtheorem{definition}{Definition}
\newtheorem{proposition}{Proposition}

\newtheorem{theorem}{Theorem}

\def\gg {\mathfrak{g}}

\def\uu {\mathfrak{u}}

\renewcommand{\(}{\begin{equation}}
\renewcommand{\)}{\end{equation}}
\newcommand{\bea}{\begin{eqnarray}}
\newcommand{\eea}{\end{eqnarray}}
\newcommand{\R}{{\mathbb R}}
\newcommand{\C}{{\mathbb C}}
\newcommand{\Z}{{\mathbb Z}}
\newcommand{\Q}{{\mathbb Q}}

\def\proof {{Proof.}\hspace{7pt}}

\def\endofproof {\hfill{$\Box$}\\}

\begin{document}

\title{ Fivebrane Structures}

\author{
  Hisham Sati,
  Urs Schreiber  and 
  Jim Stasheff} 
  \thanks{hisham.sati@yale.edu, schreiber@math.uni-hamburg.de, jds@math.upenn.edu}


\maketitle

\begin{abstract}
 We study the cohomological physics of fivebranes
 in type II and heterotic string theory. We give an interpretation of
 the one-loop term in type IIA, which involves the first and second
 Pontrjagin classes of spacetime, in terms of obstructions to 
having bundles with certain structure groups.  
 Using a generalization of the Green-Schwarz anomaly 
cancelation in heterotic string theory which demands the target space
to have a \emph{String structure},
we observe that the ``magnetic dual'' version of 
the anomaly cancelation condition
can be read as a higher analog of String structure, 
which we call {\it Fivebrane structure}. 
This involves lifts of orthogonal and unitary structures
through higher connected covers which are not just
3- but even 7-connected. 
We discuss the topological obstructions to the existence
of Fivebrane structures.
The dual version of the anomaly cancelation 
points to a relation 
of String and Fivebrane structures under electric-magnetic duality.

\end{abstract}

\newpage

\tableofcontents


\section{Introduction}

Cohomological physics  began
with Gauss in 1833, if not sooner (cf. Kirchoff's laws).
The cohomology referred to in Gauss' work was that of differential forms,
div, grad, curl and especially Stokes Theorem (the de Rham complex).
 Gauss explicitly defined the linking number of
two circles imbedded in 3-space by an integral defined
in terms of the electromagnetic
effect of a current circulating in one of the circles

\medskip
Although Maxwell's equations were for a long time expressed only
in coordinate dependent form, they were recast in a particularly
attractive way in terms of differential forms on Minkowski space.
More subtle
differential geometry and implicitly characteristic classes occurred
visibly in Dirac's magnetic monopole \cite{Dirac}, which lived in a
$U(1)$ bundle over ${\R}^3-0$.  The magnetic charge was given by the
first Chern number; for magnetic charge $1$, the monopole lived in the Hopf
bundle, introduced that same year 1931 by Hopf \cite{hopf}, though it seems to
have taken some decades for that coincidence to be recognized
\cite{greub-petry}.  Thus were characteristic classes (and by
implication the cohomology of Lie algebras and of Lie groups) introduced into
physics.

\vspace{3mm}
In the case of electromagnetic theory, the Lie group was just $U(1)$, but more general Lie groups were involved in Yang-Mills theory. There was a linguistic barrier between physicists and mathematicians 
(see Yang's reminiscences \cite{Yang}), which was breached when it was realized that the physicist's
gauge potential
 and field strength were, respectively, the mathematicians connection and curvature from differential geometry.

\vspace{3mm}
A major development occurred when Dirac's theory of the electron required -- in modern 
language-- lifting from the special orthogonal group 
$\mathrm{SO}(n)$ to the spinor group 
$\mathrm{Spin}(n)$, which corresponds to ``killing'' the first
homotopy group $\pi_1(\mathrm{SO}(n))$ of $\mathrm{SO}(n)$. 

\begin{figure}[h]
  \begin{center}
 \begin{tabular}{c|ccccccccc|l}
   $k$ && $0$ & $1$ & $2$ & $3$ & $4$ & $5$ & $6$ & $7$ &
   \\
   \hline
   &&&&&&&&&&
   \\
   $\pi_k(\mathrm{O}(n))$ && $\mathbb{Z}_2$ & $\mathbb{Z}_2$ & $0$ & $\mathbb{Z}$ & $0$ & $0$ & $0$ &
      $\mathbb{Z}$
   &
   \mbox{fundamental groups}
   \\
   \hline
   &&&&&&&&&
   \\
   $\mathrm{O}(n)\langle k\rangle$
   &
   $\mathrm{O}(n)$
   &
   $\mathrm{SO}(n)$
   &
   $\mathrm{Spin}(n)$
   &
   &
   $\mathrm{String}(n)$
   &
   &
   &
   &
   $\mathrm{Fivebrane}(n)$
   &
   \mbox{higher connected covers}
   \\
   &
   \multicolumn{9}{l}{
     \xymatrix@C=4pt{
        \hspace{10pt}
        \ar@{|->}@/_1pc/[rr]_{\mathrm{kill}\; \pi_0}&&
        \hspace{20pt}
        \ar@{|->}@/_1pc/[rrr]_{\mathrm{kill}\; \pi_1}&&&
        \hspace{14pt}
        \ar@{|->}@/_1pc/[rrrrr]_{\mathrm{kill}\; \pi_3}&&&&&
        \hspace{14pt}
        \ar@{|->}@/_1pc/[rrrrrrrrrrr]_{\mathrm{kill}\; \pi_7}
        &&&&&&&&&&&
     }}
 \end{tabular}
 \end{center}
 \caption{
   \label{homotopy groups of O(n)}
   {\bf Homotopy groups of $\mathrm{O}(n)$ and its higher connected covers} for $n > k+1$.
 }
\end{figure}

Much later it was found that
in  string theory a further step is needed, 
namely lifting the spinor group to the 
String group $\mathrm{String}(n)$, giving rise to String structures.
Such structures interpreted in terms of 
 the vanishing condition of the worldsheet anomaly of a superstring
 were first identified by Killingback \cite{Kill} 
and shortly afterwards amplified by Witten \cite{Wit}
in the context of the index theory of Dirac operators on loop space.
A nice review and differential geometric description can be found in \cite{MurrayStevenson}. 
In all these articles,
the String structure is regarded as a lift of 
an $L \mathrm{Spin}(n)$-bundle
over the free loop space $LX$ through the Kac-Moody central extension 
$\widehat {L \mathrm{Spin}}(n)$-bundle. 
It was in \cite{ST}
where it was pointed out that this lift can also be interpreted 
as a lift of the 
original $\mathrm{Spin}(n)$-bundle down on target space $X$
to a principal bundle 
for the topological group
called $\mathrm{String}(n)$. Then it was 
realized in \cite{BCSS} (see also \cite{Henriques}) that this topological
group is  in fact the realization of the nerve of a \emph{smooth} 
(as opposed to just topological) albeit categorified group: 
the String 2-group. This paves the way for a differential geometric
treatment of String structures on the target space $X$.

\vspace{3mm}
It was known early on that 2-bundles 
\cite{Bartels,BaezSchreiber}
(aka ``crossed module bundle gerbes'' \cite{Jurco})
with structure 2-group the String 2-group
have the same classification as that
of the String-bundles considered by Stolz-Teichner
(compare \cite{Jurco} and \cite{BCSS}), and
hence, rationally \cite{CP}, 
as that of the bundles on loop space originally 
considered by Killingback and by Witten.
Detailed proofs of this have recently appeared \cite{BBK} \cite{BS}. 
Note that the
String condition is closely related to the condition 
on vanishing of the integral seventh
Stiefel-Whitney class, $W_7=0$, observed in \cite{DMW} 
and studied in connection
to generalized cohomology in \cite{KS1}. In particular, 
the string condition implies the vanishing of $W_7$ as
$W_7= Sq^3(\frac{1}{2} p_1)$, where $Sq^3$ is the 
Steendrod square that raises the cohomology degree 
by three.

\vspace{3mm}
 
The word \emph{string} being well established for maps 
from 1-dimensional manifolds, 
higher dimensional analogs are referred to as \emph{branes} 
(originally, membranes). 
The ``surface" formed by an evolving 
string is called the worldsheet and, analogously, the higher-dimensional 
volumes of  evolving branes are referred to as the worldvolumes. 
 The most studied types of branes in string theory are the famous 
D-branes, which couple to the Ramond-Ramond (RR) fields. We will not
be dealing with D-branes, and hence with RR fields, in this paper. 
In addition to D-branes, there are the Neveu-Schwarz (NS) branes
that couple to the NSNS fields. The fundamental string couples to the
B-field $B_2$, whose curvature is $H_3$ ( $=dB_2 +$ non-exact, locally).
The Hodge
dual of $H_3$ in ten dimensions is $H_7$, which can be viewed as the
curvature of a degree six `potential' $B_6$, to which couples 
an extended object, called the NS 5-brane.
Thus, the 5-brane is the \emph{magnetic dual} to the string 
(in the sense familiar in String theory, which mathematicians
can find nicely described in \cite{Freed} in terms
of differential characters, which are just another
way of talking about higher line bundles with connections).
 It is to be expected that anomaly freedom
of the spinors on the fivebrane's
worldvolume
require the target to carry an analog of a spin
structure even more strict than a String structure.  
There is
a known formula for the dual to the Green-Schwarz anomaly.
This formula is discussed in section \ref{dualGS}.
The same formula can be deduced from anomaly freedom
of the worldvolume theory of the super 5-brane \cite{LT,DDP}.

\vspace{3mm}
We observe that this known formula can be read as saying that
target space $X$ needs to admit a lift of the structure group
of $T X$ from $\mathrm{SO}(n)$ through $\mathrm{Spin}(n)$,
through $\mathrm{String}(n)$ and then further to the 7-connected
cover of $\mathrm{SO}(n)$ which we dub $\mathrm{Fivebrane}(n)$.
We define a {\it Fivebrane structure},
which is obtained 
by killing all up to and including the \emph{seventh} 
homotopy group of the 
orthogonal group. This is entirely analogous 
to the process of killing the third 
homotopy group in going from the Spin group to the String group. 
Generally, one could wonder what the even higher connected
covers of $\mathrm{SO}(n)$ would correspond to in terms of
higher brane physics.
While the notions of higher lifts are of 
course known in homotopy theory, we here discuss
the topological structure associated with 
the physics of fivebranes, much the same as in the known String case.

 \begin{definition}
 An $n$-dimensional manifold $X$ has a
 \emph{ fivebrane structure} 
 if the classifying
 map $X \to (BO)(n)$ of the tangent bundle $T X$ lifts to 
 $B{\rm Fivebrane} :=(B\mathrm{O})\langle 9 \rangle$.
$$
  \raisebox{20pt}{
  \xymatrix{
     &
     \mathrm{B O} (n) \langle 9 \rangle 
     \ar[d]
     \\
     M
     \ar[r]_{\hspace{-3mm}f}
     \ar@{..>}[ur]^{\hat f}
     &
     \mathrm{B O}(n)~.
  }
  }
$$
\label{maindef}
\end{definition}

One point we make is that the fivebrane -- as opposed to branes of other dimensions-- 
is distinguished in 10-dimensional spacetime since it does indeed lead to a structure 
that naturally generalizes 
the string structure, due to the existence of the field $H_7$, the dual field
to $H_3$. The importance of $H_7$ will be highlighted in section \ref{dualGS}.
Notice that the existence of this lift implies lifts through
all the lower connected covers, which says that 
for $X$ to have a $\mathrm{Fivebrane}$ structure it must also have
an orientation, a Spin structure and a String structure.

\vspace{3mm}
Our aim is to
\begin{itemize}
  \item
    understand the {\bf topological nature} of such 
    \emph{fivebrane structures}, i.e. identify the 
    relevant bundles and the corresponding 
    characteristic classes;
  \item
    understand their  {\bf differential geometric nature} in terms
    of characteristic \emph{forms} of such bundles (e.g. generalized connections)
    analogous to what is done for String theory.
\end{itemize}

This paper will focus more on the first, while a sequel \cite{SSS2}
will emphasize the second.
We shall show the following:
\begin{theorem}
The obstruction to lifting a String structure on $X$
    to a Fivebrane structure on $X$ is the fractional
    second Pontrjagin class
    $ \frac{1}{6}p_2(T X)$.
\end{theorem}
    This is our Proposition \ref{p26}.
    \begin{itemize}
     \item
       Here the fractional coefficient $1/6$ is the crucial 
       subtle point.
       It is explained in \ref{String structures and 
       the second Pontrjagin class}.
      \item   
        Inequivalent $\mathrm{Fivebrane}$ structures on $X$ are
        classified by a quotient of $H^7(X,\mathbb{Z})$.
        This is proposition \ref{the set of lifts}.
    \end{itemize}

But this is not the full story yet. String structures
are in general carried not just by a manifold $X$, 
but by a complex vector bundle $E \to X$
(e.g. the gauge bundle on the target space of the 
heterotic string): 
this vector bundle
has $\mathrm{String}$ structure if its second Chern class
\footnote{We review characteristic classes and characters in Appendix A.} 
cancels the fractional first Pontrjagin class of the
Spin bundle $X$: 
\bea
  \mbox{$E\to X$ has Spin structure}
  \;\;
 &\Longleftrightarrow &
  \;\;
  \left\lbrace
  \begin{tabular}{l}
    \mbox{$E$ is orientable}
    \\
    \mbox{and}
    \\
    $w_2(E) = 0$
   \end{tabular}
   \right.
  \nonumber\\
  \mbox{$E\to X$ has String structure}
  \;\;
  &\Longleftrightarrow &
  \;\;
  \left\lbrace
  \begin{tabular}{l}
    \mbox{$E$ has Spin structure} 
    \\
    \mbox{and}
    \\
    $\frac{1}{2}p_1(T X) = \mathrm{ch}_2(E)$
   \end{tabular}
   \right.
  \,.
\nonumber
\eea
If the second Chern character of $E$ vanishes, this reduces
to saying that $X$ itself has String structure.
(The fractional first Pontrjagin class is reviewed in section
\ref{spin structures and first pontrjagin}.)

\vspace{3mm}
This situation generalizes. A Fivebrane structure
should be assigned, more generally, to a vector bundle $E \to X$.
The dual Green-Schwarz mechanism indicates that the 
condition on $E$ to have a Fivebrane structure is essentially
that the fourth Chern class of $E$ cancels the second Pontrjagin
class of $X$, but there are corrections by decomposable classes
(because ${\rm ch}_2=c_2 + \mbox{decomposable classes}$):
$$
  \hspace{-.5cm}
  \mbox{$E\to X$ has Fivebrane structure}
  \;\;
  \Longleftrightarrow
  \;\;
  \left\lbrace
  \begin{tabular}{l}
    \mbox{$E$ has String structure} 
    \\
    \mbox{and}
    \\
    $\frac{1}{6}p_2(T X) = \mathrm{ch}_4(E)$
   + \mbox{decomposable classes}
   \end{tabular}
   \right.
  \,.
$$

\begin{table}[h]
  \begin{center}
  \begin{tabular}{lr|ll}
    \multicolumn{2}{c}{\bf higher spin-like structure}
    &
    \multicolumn{2}{c}{ \bf defining condition}
    \\
    \hline
    &
    \\
    $\mathrm{Spin}$ structure on manifold & $X$ & 
     orientation and & $w_2(T X) = 0$ 
    \\
    \hline
    &
    \\
    $\mathrm{Spin}$ structure on gauge bundle & $E \to X$
    &
    orientation and & $w_2(E) = 0$
    \\
    \hline
    &
    \\
    $\mathrm{String}$ structure on manifold & $X$ & 
     $\mathrm{Spin}$ structure and & $\frac{1}{2}p_1(T X) = 0$ 
    \\
    \hline
    &&
    \\
    $\mathrm{String}$ structure on gauge bundle & $E \to X$ & 
    $\mathrm{Spin}$ structure and &
     $\frac{1}{2}p_1(T X) = \mathrm{ch}_2(E)$ 
    \\
    \hline
    &&
    \\
    $\mathrm{Fivebrane}$ structure on manifold & $X$ & 
     $\mathrm{String}$ structure and
     &
     $\frac{1}{6}p_2(T X) = 0$ 
    \\
    \hline
    &&
    \\
    $\mathrm{Fivebrane}$ structure on gauge bundle & $E \to X$ & 
    $\mathrm{String}$ structure and
     &
     $\frac{1}{6}p_2(T X) = 8 \mathrm{ch}_4(E) + \mathrm{decomposables}$ 
 \\
 \hline
  \end{tabular}
  \end{center}
  \vspace{2mm}
  \caption{
  {\bf Higher Spin-like structures} on manifolds and on 
  ``gauge bundles'' over a manifold $X$. Here ``gauge bundle'' 
  technicaly just means: complex vector bundle. 
  The first four entries are
  well established; we are concerned here with the last two 
  entries.
  }
\end{table}

We will study fivebrane structures arising from anomaly 
polynomials in two theories: in 
type IIA string theory 
and in heterotic string theory. In the first, we will have
a fivebrane structure for the tangent bundle of the space
itself, since there is no gauge bundle. We address this
in full detail. In the second theory, we have, in addition, 
a gauge bundle. The corresponding condition will involve 
both the tangent bundle and this gauge bundle. We 
are still able to give a description when the two 
separate conditions are imposed: that {\it separately}
the tangent bundle and the gauge bundle admit 
a fivebrane structure. In the general case we encounter
a difficulty, namely that the factors -- the second 
Pontrjagin classes-- come with different coefficients. 
We explain this difficulty and provide (partial) solutions in 
section \ref{congruence} and towards the end of 
 section \ref{ineq}. We will be concerned in this paper primarily 
 with the cohomological aspects
of Fivebrane structures rather than the more subtle aspects related 
to K-theory which we hope to address elsewhere.

\vspace{3mm}
The presence of the decomposable cohomology classes 
(products of nontrivial
cohomology classes) appearing on the right can be deduced 
from the dual Green-Schwarz formula. Notice that
a decomposable characteristic class necessarily suspends
to a trivial group cocycle. So despite the appearance of
those decomposable classes, Fivebrane structures are
essentially controlled just by the 7th Lie algebra cohomology
of $\mathfrak{so}(n)$, just as String structures are controlled
by the third Lie algebra cohomology (compare \cite{SSS1}).
We shall discuss in \cite{SSS2} that the decomposable 
characteristic classes in the dual Green-Schwarz formula
affect the differential class (the connection) but not the 
integral class of the smooth Chern-Simons 7-bundle 
whose non-triviality obstructs the existence of 
$\mathrm{Fivebrane}$ structures. 
The description in terms
of anomalies will be given in section 
\ref{anomaly cancellation in string theory}.

\vspace{3mm}
The notion of and notation for connected covers is explained towards the
end of section \ref{The homotopy groups of SO}.

\vspace{3mm}
\noindent {\bf
Outlook: differential geometric interpretation.}
Just as String 2-bundles with connection provide a differential
geometric realization of  $\mathrm{String}$ structure, the higher
analog of a $\mathrm{Spin}$ bundle with connection, there are
analogous higher bundles with connection
giving a differential geometric realization of
$\mathrm{Fivebrane}$ structures.
Following the counting pattern these deserve to be called 
$\mathrm{Fivebrane}$-6-bundles or $\mathrm{Fivebrane}$-5-gerbes
with connection, but the reader might want to think of them
just as a nonabelian version of differential cocycles with top
degree curvature form a 7-form. 
As will be described in \cite{SSS2}, these nonabelian differential
cocycles can be obtained from integrating the $L_\infty$-algebra
connections of \cite{SSS1}. There the Lie 2-algebra 
$\mathfrak{string}(n)$ 
(an $L_\infty$-algebra concentrated in the lowest two degrees)
governing String bundles with connection
had been accompanied by the Lie 6-algebra 
$\mathfrak{fivebrane}(n)$ (concentrated in the lowest six degrees)
which plays the corresponding role for Fivebrane structures.

\vspace{3mm}
In the sense of integration of Lie algebras to Lie groups, the
$\mathfrak{string}(n)$ Lie 2-algebra had been integrated 
in \cite{Henriques,BCSS} to the
$\mathrm{String}(n)$ Lie 2-group, the structure Lie 2-group
of String 2-bundles. This integration procedure can be thought 
of as follows: 
for any $L_\infty$-algebra $\gg$ there is a notion of 
$\gg$-valued differential forms \cite{SSS1} and there is a generalized 
smooth space, $S(\mathrm{CE}(\gg))$, which is the \emph{classifying}
space for $\gg$-valued forms in that smooth maps from any other
smooth space $Y$ into it are
in bijection with $\gg$-valued forms on $Y$:
\(
  \mathrm{Hom}_{\mathrm{SmoothSpaces}}(Y,S(\mathrm{CE}(\gg)))
  \simeq
  \Omega^\bullet(Y,\gg)
  :=
  \mathrm{Hom}_{\mathrm{DGCA}s}(\mathrm{CE}(\gg),\Omega^\bullet(Y))  
  \,.
\) 
The $\infty$-groups integrating given $L_\infty$-algebras $\gg$
considered in \cite{Getzler,Henriques} are 
obtained essentially by forming the
$\infty$-path groupoid in the sense of Kan complexes of this
classifying space, which is nothing but the
singular simplicial complex of $S(\mathrm{CE}(\gg))$. This 
corresponds to a ``maximally weakened'' integration. 
Alternatively,
one can form the strict globular (instead of simplicial) 
$n$-path groupoid $\Pi_n$ (whose $n$-morphisms are
homotopy classes of $n$-paths and whose $(k < n)$-morphisms are
\emph{thin} homotopy classes of $k$-paths) which yields the ``maximally strictified''
$n$-group integrating the given Lie $n$-algebra. This way 
one obtains a strict version of the String Lie 2-group
\(
  \mathbf{B}\mathrm{String}(n)
  :=
  \Pi_2\left( S(\mathrm{CE}(\mathfrak{string}(n))) \right)
\)
and the strict version of the Fivebrane Lie 6-group
\(
  \mathbf{B}\mathrm{Fivebrane}(n)
  :=
  \Pi_6(S(\mathrm{CE}(\mathfrak{fivebrane}(n))))
  \,.
\)
Here on the left the notation $\mathbf{B}G$ indicates 
a (strict) $n$-groupoid with a single object corresponding
to the (strict) $n$-group $G$. The notation is such that
taking geometric realization of nerves $|\cdot|$ to produce
topological spaces we have
\(
  |\mathbf{B}G| = B |G|
  \,,
\)
where on the right $|G|$ is a topological group and $B |G|$ its
classifying space \cite{BS}.
The strict realizations of the String and Fivebrane
Lie $n$-groups allow us
to employ Ross Street's 
general theory of descent \cite{StreetDesc} to
formulate String 2- and Fivebrane 6-bundles with connection
in terms of nonabelian differential cocycles \cite{ndclecture}.

\section{The context}

Before getting to our main discussion, we first 
indicate the general context in which these questions
arise.
The word \emph{string} being well established for maps 
from 1-dimensional manifolds, 
higher dimensional analogs are referred to as \emph{branes} 
(originally, membranes). 
The ``surface" formed by an evolving 
string is called the worldsheet and, analogously, the higher-dimensional 
volumes of  evolving branes are referred to as the worldvolumes. 
We have an ``$n$-particle", otherwise known as an $(n-1)$-brane, whose worldvolume
is an 
$n$-dimensional manifold $\Sigma$, or rather the image of that $n$-dimensional 
manifold under a sufficiently well behaved map
\(
  \phi : \Sigma \to X
  \label{phi}
\)
 into a target space (usually to be thought of as physical spacetime) $X$. 
On that spacetime, we have a (generalized, higher) 
bundle (thought of as an $n$-bundle or $(n-1)$-gerbe)
with (generalized) connection.  
Henceforth, we will drop the `generalized' and the $n$ unless crucial. Thinking of an ordinary bundle with connection should be sufficient. This bundle with connection encodes the data specifying 
the ``background field'' to which
that  $(n-1)$-brane ``couples''. 
Just as with ordinary connections, 
connections on an $n$-bundle can be expressed either
in terms of local $(p \leq n)$-forms on the base manifold 
or in terms of  global $(p \leq n)$-forms
on the total space \cite{SSS1}.
The local representatives of these forms are
often referred to as a ``background field'', though technically
speaking the background field is that differential form
datum together with the descent/gluing data that makes it
a differential cocycle. 
In particular,  we have the 3-form curvature $H_3$
of the Kalb-Ramond field in string theory,
which is, in this sense, the curvature of a 
``background field'', the ``$B$-field'', of a string theory.
Likewise, we have
the 7-form $H_7$ which is the curvature of 
a background field for the fivebrane.

\subsection{$\Sigma$-models}

We are concerned with the mathematical structure which
is supposed to model the physics of \emph{charged $n$-particles},
usually known as \emph{charged $(n-1)$-branes} or as
\emph{quantum field theories of $\Sigma$-model type}.
Such a $\Sigma$-model is specified by choosing
\begin{itemize}
  \item
    a ``space'' $X$, called the 
    \emph{target space};
  \item
    a ``space'' (or class of such) $\Sigma$, called
     the \emph{parameter space} or called the \emph{worldvolume};
  \item
    the mapping space $\mathrm{Maps}(\Sigma,X)$ called the
    \emph{space of fields} or the \emph{configuration space}
    or sometimes the \emph{moduli space} (the latter is usually a quotient);
  \item
    on the target space a \emph{differential $n$-cocycle}
    $\nabla$, i.e. a higher generalization of a fiber bundle
    with connection, called the \emph{background field};
  \item
    a prescription for how to interpret the push-forward 
    of the the pullback $\mathrm{ev}^* \nabla$ 
    along the projection $pr_1$ onto $\Sigma$ in the
    correspondence diagram
    \(
      \label{the correspondence}
      \raisebox{20pt}{
      \xymatrix{
        & 
        \Sigma \times \mathrm{Maps}(\Sigma,X)
        \ar[dl]_{pr_1}
        \ar[dr]^{\mathrm{ev}}
        \\
        \Sigma
        &&
        X
      }
      }
    \)
    called the \emph{path integral} or the 
    \emph{quantization of the $\Sigma$-model}.
\end{itemize}

\vspace{3mm}
When the parameter space $\Sigma$ is $n$-dimensional,
one thinks of this data as encoding the physics of 
$n$-fold higher analogs of particles, ``$n$-particles'', 
that propagate on $X$. The field configuration (\ref{phi}) is thought of 
as the \emph{trajectory} of such
an $n$-particle in $X$.

\begin{table}[h]
  \begin{center}
  \begin{tabular}{c|c|c|c}
     \begin{tabular}{c}
       {\bf fundamental}
       \\
       {\bf object}
     \end{tabular}  
     &
     \begin{tabular}{c}
       {\bf dimension of}
       \\
       {\bf worldvolume $\Sigma$}
     \end{tabular}
     &
     \begin{tabular}{c}
       {\bf background}
       \\
       {\bf field}
     \end{tabular}     
     &
     \begin{tabular}{c}
       {\bf curvature}/\\
       {\bf field strength}
     \end{tabular}          
     \\
     \hline
     &&&
     \\
     $n$-particle & $n$&  $n$-bundle
     &
     $(n+1)$-form
     \\
     &&&
     \\
     $(n-1)$-brane &  $n$&  $(n-1)$-gerbe
     &
     $(n+1)$-form
\\
\hline
  \end{tabular}
  \end{center}
  \vspace{2mm}
  \caption{
    {\bf The two schools of counting} higher dimensional 
    structures. Here $n = \mathrm{dim}(\Sigma)$ is in $\mathbb{N} = \{1,2,\cdots\}$.
   The $n=1$ case: A particle is a 0-dimensional spatial object moving in time, so 
that the `worldvolume' $\Sigma$ is $(0+1)$-dimensional and hence 
1-dimensional in our terminology. The corresponding background field 
is encoded in a 1-bundle or 0-gerbe with curvature 2-form. The generalization to
strings ($n=2$) and higher extended objects, i.e. the $n$-particles, 
works analogously. One can also include the case $n=0$ corresponding to 
``instantons", although the corresponding geometric description will be different. 
  }
\end{table}

One says that the $n$-particle \emph{couples} to the 
background field $\nabla$ or that it is \emph{charged under}
the background field. The terminology is entirely motivated from
the familiar case of ordinary electromagnetically 
charged $1-$particles: the electromagnetic background field
$\nabla$ which they couple to is modeled by a vector bundle 
(a line bundle in this case) with connection.
For $n=2$ one speaks of 
``strings''. String theory strictly is the study of 
those $n=2$ $\Sigma$-models with a special restriction 
for what the ``path integral''
is allowed to be. Technically, string theory
 is required to encode
a 2-dimensional superconformal field theory of central charge
15.  This condition, however, is of no real relevance for
our discussion here, which pertains to all $\Sigma$-models
which generalize the ``spinning 1- particle''. 

\vspace{3mm}
Some of the most interesting ideas concerning such $\Sigma$-models
have originally been thought by Dan Freed:
The interpretation of background fields and of charges
as differential cocycles is nicely described and worked
out in \cite{Freed,HS}, where the mathematically inclined
reader can find rigorous interpretations,
in terms of differential cohomology, of the \emph{abelian} kinds
of ``background fields'' and related ``anomalies''
in string theory with which we are concerned here,
which include the Chern-Simons 3- and 7-bundles with 
connection obstructing the $\mathrm{String}(n)$ and
$\mathrm{Fivebrane}(n)$ lifts, but not these lifts 
themselves. (For those,
one would need \emph{nonabelian}
differential cohomology \cite{SSS2}).

\vspace{3mm}
\noindent {\bf Integration as pushforward.}
The interpretation of quantization and of the path integral
as an operation on higher categorical structures has
first been explored in \cite{FreedI,FreedII}. Integration 
as a push-forward operation plays a promint role
in recent developments by Stolz and Teichner
and by Hopkins et al. 
Let us take for instance the simple
toy example case where the background field 
$\nabla$ is a vector bundle (without connection) and where
$\Sigma$ is a point, i.e. $n=0$. In this case $pr_1$ is the map from 
${\rm pt} \times {\rm Maps}({\rm pt}, X)=X$ to the point. 
Integration over the fiber in this case is just integration over
$X$ itself, and
the ordinary push-forward gives
the space of sections of the original vector bundle, as 
the point is varied over $X$. 
That
reproduces indeed the desired ``quantization over the
point'' and can, following \cite{FreedI,FreedII}, 
be regarded as the codimension 1 part of the full path integral
for $n=1$.
Stolz and Teichner describe a 
variation of this which involves push-forward of 
K-theory classes to the point, which then classifies 
connected components of all (supersymmetric) 1-dimensional
$\Sigma$-models.
This shows that, while a fully satisfactory 
mathematical interpretation of the quantization of 
$\Sigma$-models is to date still 
an open question, a coherent picture, revolving around
the correspondence (\ref{the correspondence}),
is beginning to emerge. The ``higher spin-like structures''
on target space $X$ discussed here are believed to ensure
the existence of the quantization step in the case that
the $\Sigma$-model generalizes that describing spinning 
1-particles.

\vspace{3mm}
The sigma model for the string (i.e. $\Sigma=$ string worldsheet) 
has been studied extensively in the literature and provides a consistent
model both at the classical
and the quantum levels. Due to string/fivebrane duality in ten dimensions
\cite{Str} \cite{CHS} \cite{DL1} \cite{DL2},
one expects that there should be a formulation of string theory via 
fivebrane sigma models (i.e. $\Sigma=$ fivebrane worldvolume). 
While this is expected, we point out that such a program has not
been fully completed. However, there are works that point in that
direction \cite{DDP}. At the least for the gravitational 
parts-- i.e. without the gauge bundle-- for the fivebrane,
there are models
in which the anomalies from the worldvolume theory match
the expression of the polynomials in the Pontrjagin classes
\cite{DDP} \cite{LT}.
That is enough for our purposes since our main focus is
the
cohomological
structure resulting from lifts of the tangent bundle of spacetime
and, after all, it is not our aim to write down a full quantum 
fivebrane sigma model action.
  The paper \cite{LT} also highlights some of the 
difficulties encountered, but also gives partial resolutions.  

\vspace{3mm}
The reader not further concerned with string theoretic reasoning
might proceed to section \ref{fivebrane structures}.

\subsection{Background fields}

Independently of how the ``background field'' $\nabla$ is
modeled, it should locally be encoded by differential 
form data. See table \ref{fields table}

\begin{table}[h]
  \begin{center}
   \small
  \begin{tabular}{llll}
    {\bf $n$-particle}
    &
    {\bf background field}
    &
    {\bf global model}   
    &
    {\bf \begin{tabular}{l}local differential \\form data\end{tabular}}
    \\
\hline
    $(1-)$particle & electromagnetic field 
    &
    \begin{tabular}{l}
      line bundle with connection/\\
      Cheeger-Simons differential 2-character/\\
      Deligne 2-cocycle
    \end{tabular}
    &
    \begin{tabular}{l}
      connection 1-form:\\ $A \in \Omega^1(Y)$
      \\
      curvature 2-form:\\ $F_2 := dA \in \Omega^2_{\mathrm{closed}}(Y)$
    \end{tabular}
    \\
    \hline
    \\
    \begin{tabular}{c}
       string 
       \\
       (2-particle)
       \\
       (1-brane)
    \end{tabular}
    & Kalb-Ramond field 
    &
    \begin{tabular}{l}
      line 2-bundle with connection/\\
      bundle gerbe with connection (``and curving'')/\\
      Cheeger-Simons differential 3-character/\\
      Deligne 3-cocycle
    \end{tabular}
    &
    \begin{tabular}{l}
      connection 2-form: \\ $B \in \Omega^2(Y)$
      \\
      curvature 3-form: \\ $H_3 := dB \in \Omega^3_{\mathrm{closed}}(Y)$
    \end{tabular}
    \\
    \hline
    \\
    \begin{tabular}{c}
       membrane
       \\
       (3-particle)
       \\
       (2-brane)
    \end{tabular}
    & 
     \begin{tabular}{c}supergravity \\ 3-form field\end{tabular}
    &
    \begin{tabular}{l}
      line 3-bundle with connection/\\
      bundle 2-gerbe with connection (``and curving'')/\\
      Cheeger-Simons differential 4-character/\\
      Deligne 4-cocycle
    \end{tabular}
    &
    \begin{tabular}{l}
      connection 3-form: \\ $C \in \Omega^3(Y)$
      \\
      curvature 4-form: \\ $G_4 := dC \in \Omega^4_{\mathrm{closed}}(Y)$
    \end{tabular}
    \\
    \hline
    \\
  \end{tabular}
  \end{center}
  \caption{
    {\bf Simple (abelian) examples 
     for $n$-particles and the background fields
     they couple to}.
    The background fields are often addressed in terms of 
    the symbols used for their local form data:
    the Kalb-Ramond field is known as the ``$B$-field''
    with its ``$H_3$ field strength'' . Similarly 
    one speaks of the ``$C_3$-field'' and its field strength
    ``$G_4$'', etc. This reflects
    the historical development, where the local differential
    form data was discovered first and its global
    interpretation only much later. (See also the 
    remark on anomalies at the beginning of 
    section \ref{dualGS}).
  }
  \label{fields table}
\end{table}

\vspace{3mm}
All the relevant background fields that have
been considered are locally controlled by some 
$L_\infty$-algebra $\gg$, and the local differential
form data can always be considered as encoding 
differential forms $A \in \Omega^\bullet(Y,\gg)$ with 
values in the Lie algebra $\gg$ \cite{SSS1}.
In the case of abelian differential cocycles, these
$L_\infty$-algebras are all of the form $b^{n-1}\uu(1)$:
the higher dimensional versions of $\uu(1)$.



\subsection{Charges}

Just as an ordinary 1-bundle may be trivialized by a \emph{non-trivial} section,
which one may think of as a
 ``twisted 0-bundle'', higher
$n$-bundles may be trivialized by ``higher sections'' which
are called ``twisted $(n-1)$-bundles''.
One says the twisted $(n-1)$-bundle is ``twisted by''
the corresponding $n$-bundle.
A beautiful description of this situation for abelian
$n$-bundles with connection in terms of differential
characters is given in \cite{Freed,HS}.
Twisted nonabelian 1-bundles have been studied in detail
under the term ``bundle gerbe modules'' 
\cite{BCMMS}.
Twisted non-abelian 2-bundles have first been considered in
\cite{AschieriJurco,Jurco} under the name
``twisted crossed module bundle gerbes''. 
In terms of the $L_\infty$-connections
considered in \cite{SSS1} twisted $n$-bundles with connections
are the connections for $L_\infty$-algebras arising as
mapping cone $L_\infty$-algebras $(b^{n-1}\uu(1) \to \hat \gg)$.

\vspace{3mm}
By comparing the formalism here with the situation of
ordinary electromagnetism, one can regard the twisting 
$n$-bundle as encoding the presence of ``magnetic charge''.
This, too, is nicely explained at the beginning of \cite{Freed}.
Accordingly, where an untwisted $(n-1)$-bundle has a curvature
$n$-form $H_{n}$ which is closed, a twisted
$(n-1)$-bundle has a curvature $n$-form which is 
``twisted by'' the curvature $(n+1)$-form $G_{(n+1)}$ of the 
twisting $n$-bundle
\(
  d H_{n} = G_{n+1}
  \,.
\)
Indeed, for a twisted $(n-1)$-bundle the curvature
$H_n$
is locally no longer the differential of the connection,
$d B_{n-1} = H_n$, but receives a contribution from the
connection $n$-form $B_n$ of the twisting $n$-bundle
\(
  H_n = d B_{n-1} + B_n
  \,.
\)
The archetypical example is that of ordinary magnetic
charge: as Maxwell discovered in the 19th century, 
in the presence of 
magnetic charge, which in four dimensions
is modelled by a 3-form $H_{3} = \star j_1$, the 
electric field strength 2-form $F_2$ is no longer closed
\(
  d F_2 = H_3
  \,.
\)
When Dirac later discovered at the beginning of the
20th century that $H_3$ has to have
integral periods (``quantization of magentic charge''), 
the first 2-categorical structure in
physics had appeared: the magnetic 
torsion 2-bundle / bundle-gerbe
with deRham class $H_3$. It seems that this was first
explicitly realized in \cite{Freed}.

\vspace{3mm}
The next example of this kind received such a  
great amount of attention
that it came to be known as the initiation of the 
``first superstring revolution'': the Green-Schwarz 
anomaly cancelation mechanism \cite{GS}.
The interaction of gauge theory with type I supergravity theory leads
to the (low energy limit of the) heterotic theory which has a rich mathematical 
structure. The 
 \emph{Chapline-Manton coupling} \cite{CM} 
 amounts essentially to equating the curvature
 $H_3$ of the $B$-field in type I  with the 
 Chern-Simons
 3-form of the connection on the gauge bundle. More precisely,
 in terms of the the virtual difference of two Chern-Simons
3-bundles (Chern-Simons 2-gerbes), locally, one has
\(
  d H_3 = d \mathrm{CS}(\omega) - d \mathrm{CS}(A)
\label{csdiff}
\)
for $\omega$ and $A$ the local connection 1-forms of 
a Spin and complex vector bundle, respectively and
$\mathrm{CS}(-)$ denoting the corresponding Chern-Simons
3-forms. 
 Thus this leads to nontrivial and rich structures both physically and 
 mathematically. See for instance \cite{Freed} \cite{BCRS} and \cite{Clin} for
 geometric treatments.

\vspace{3mm}
As nicely explained by \cite{Freed},
in the ``higher gauge theory'' given by the
effective supergravity target space theory of the heterotic string,
the supergravity $C$-field with curvature 4-form 
\footnote{In all of this paragraph, what we mean by the $C$-field
and its corresponding curvature $G_4$ is their topological part,
i.e. the shift of $G_4$ in the formula 
$[G_4] -\frac{1}{4}p_1= a \in H^4(Y^{11}, \Z)$ 
\cite{Flux}, and not the fields themselves --the `$p_1$ part' 
and not the `$G_4$ part' of $a$.}
$G_4$ had to be ``trivialized''
by the Kalb-Ramond field with curvature 3-form $H_3$,
or conversely the Kalb-Ramond field had to be ``twisted''
by the supergravity curvature 4-form  $d H_3 = G_4$
Moreover, $G_4$ has to be the curvature of the virtual difference 
implicit in (\ref{csdiff}). 
The Green-Schwarz anomaly cancelation condition can
hence be read, equivalently, as saying that
\begin{itemize}
  \item
    the supergravity $C$-field trivializes over the 
    10-dimensional target of the heterotic string;
  \item
    $G_4$ is the \emph{magnetic 5-brane charge} 
    which the electric heterotic string couples to;
  \item
   the Kalb-Ramond field is \emph{twisted} by the supergravity  
   $C$-field.
\end{itemize}

\vspace{3mm}
This anomaly cancelation has
an interesting description
in terms of
a string structure \cite{Kill} \cite{CP}. 
This can also be interpreted as a spin structure
on the (free)
loop space of 
the ten-dimensional 
spacetime \cite{Wit} \cite{Se} \cite{ST}. 
This free loop space can be 
viewed as the configuration space of the string. From a topological point 
of view, the string structure is equivalent to lifting the structure group of the tangent 
bundle of spacetime from ${\rm Spin}(10)$ to its three-connected 
cover ${\rm String}(10)$, obtained by killing the first three homotopy 
groups. The latter is infinite-dimensional
but can be captured by some
finite-dimensional constructions \cite{ST} \cite{BCSS}.

\vspace{3mm}
There is no particular reason to prefer ``electric charge''
over ``magnetic charge'': in the presence of a Riemannian
structure,
the Hodge star dual of an ``electric'' 
field strength $H_{n+1}$ may be interpreted as a field
strength itself, in which case it is called the
``magnetic field strength'' $H_{d-n-1} := \star H_{{n+1}}$.
Just as the original field strength $H_n$ coupled to
an ``electric'' $n$-particle,
the dual field strength couples to a 
``magnetic'' $(d-n-2)$-particle. 
Such electric-magnetic duality is at the heart of what
is known as ``S-duality'' for super Yang-Mills theory, 
which has recently been argued 
\cite{KapustinWitten} to be the heart of 
geometric Langlands duality.
It is only for electric 1-particles in $d=4$ dimensions that
their magnetic dual is again a 1-particle. The magnetic
dual of the 2-particle in 10 dimensions is the 6-particle.
In other words: the magnetic dual of the string is the 
5-brane.

\vspace{3mm}
Here our starting point is to look at the above situation with
the electric string in the presence of magnetic 5-brane charge in the
dual formulation, where the magnetic 5-brane couples to 
a 6-form field with field strength $H_7$ in the presence of
electric string charge, which is then given by an 8-form $I_8$.

\begin{table}[h]
  \begin{tabular}{cccc|ccccccc}
   \hline
    &
    {\begin{tabular}{c}
      electric
      \\
      field strength 
      \\
      coupled to
      \\
      fundamental 
      \\
       electric
      \\
      $(n-1)$-brane
    \end{tabular}
    }
    &&
    {\begin{tabular}{c}
      magnetic 
      \\
      $(d-n-3)$-brane  
      \\
      current
    \end{tabular}
    }  
    &&
    {\begin{tabular}{c}
      magnetic
      \\
      field strength 
      \\
      coupled to
      \\
      fundamental 
      \\
       magnetic
      \\
      $(d-n-3)$-brane
    \end{tabular}
    }
    &&
    {\begin{tabular}{c}
      electric
      \\
      $(n-1)$-brane  
      \\
      current
    \end{tabular}
    }  
    \\
    \hline 
    &&&&&&
    \\
    $d$
    &
    $F_{n+1}$
    &$=$&
    $j_B$
    &
    $d$
    &
    $\star F_{n+1}$
    &
    $=$
    &
    $j_E$
    \\
    \hline
    \hline
    &&&&&
    \\
    &&&&&
    \\
    &
    {\begin{tabular}{c}
      electric
      \\
      KR field
      \\
      coupled to
      \\
      fundamental 
      \\
      string
    \end{tabular}
    }
    &&
    {\begin{tabular}{c}
      magnetic 
      \\
      $5$-brane  
      \\
      current
    \end{tabular}
    }  
    &&
    {\begin{tabular}{c}
      magnetic
      \\
      dual KR field
      \\
      coupled to
      \\
      fundamental 
      \\
      $5$-brane
    \end{tabular}
    }
    &&
    {\begin{tabular}{c}
      electric
      \\   
      string
      \\
      current
    \end{tabular}
    }  
    \\
    \hline 
    &&&&&&
    \\
    $d$
    &
    $H_3$
    &$=$&
    $\underbrace{\frac{1}{2}p_1(\omega) - \mathrm{ch}_2(A)}_{= I_4}$
    &
    $d$
    &
    $H_7$
    &
    $=$
    &
    \raisebox{-12pt}{
    \begin{tabular}{l}
      $\frac{1}{48}p_2(\omega) - \mathrm{ch}_4(A)$
      \\
      $\underbrace{ + \mathrm{decomposables}}_{I_8}$
     \end{tabular}
     }
     \\
     \hline
  \end{tabular}
  \vspace{2mm}
  \caption{
    The Green-Schwarz formula and its dual version, 
    with its interpretation in terms of electric
    strings and their magnetic 5-brane duals.
    The electric current $I_4$ and the magnetic current
    $I_8$ are both fixed such that the anomaly they
    produce cancels the anomaly from the fermions in the
    theory (see section \ref{anomaly cancellation in string theory}).
  }
\end{table}

\vspace{3mm}
Type I supergravity admits a formulation in terms of the potential $B_6$
corresponding to the field $H_7$, which is Hodge dual to $H_3$ in ten 
dimensions \cite{Cham}. There is also a corresponding anomaly cancelation 
procedure for this dual theory which makes use of a degree seven analog
of the Chapline-Manton coupling \cite{Gates} \cite{SS}. This process is 
also mathematically rich and has been treated in \cite{Freed} from 
a K-theoretic point of view and, in fact, in a duality-symmetric fashion,
i.e. including both fields on equal footing and at the same time. 
The magnetic dual discussion of the Green-Schwarz mechanism 
\cite{Gates} \cite{SS} leads us to consideration of
a twisted 6-bundle with field strength $H_7 = \star H_3$,
which is twisted by a certain 7-bundle whose
field strength eight-form is a sum of two higher characteristic 
classes plus some mixed terms -- see equation (\ref{dH7}). 
This is the formula which we shall refer to as the
\emph{dual Green-Schwarz anomaly cancellation condition}
and take as the starting point of our discussion.

\section{The Dual Green-Schwarz Mechanism and Higher Chern-Simons Forms}

\label{dualGS}

Here we recall the string theoretic results
which indicate which characteristic classes
of a manifold $X$ have to vanish in order for the
manifold to qualify as the target for the propagation of
a 5-brane. 
For the main mathematical point that we make in section
\ref{fivebrane structures},
  it is
sufficient to note here that there is motivation, from formal
high energy physics,
for studying the condition (\ref{dH7}), below, on characteristic
classes of complex vector bundles over $X$. 
The reader not further concerned with string theoretic reasoning
might just want to note this equation and the observations following
it and then proceed to section \ref{fivebrane structures}.

\subsection{Anomaly cancelation in string theory}
\label{anomaly cancellation in string theory}

There are several somewhat different phenomena which are
called \emph{anomalies} in physics, but they usually
all refer to issues of global topological twists.
Physicists are used to develop their concepts in terms
of local data and many times implicitly assume that this is sufficient.
Generically, the anomaly in our context refers to an inconsistency in 
the topological assumptions taken for the underlying space or for 
bundles on that space.
For instance, if one accepts that spinors are sections of
Spin bundles, then it is obvious that their existence 
requires the underlying manifold to have a Spin structure.
But one way to discover this from the point of view of physics is,
as nicely described in \cite{WittenOld},
to start with a naive action functional for a spinning particle
and then to discover that it is ill defined globally unless the
target has a Spin structure.
Entirely analogous considerations lead to String structures
as the ``anomaly cancelation conditions'' for superstrings,
known as the Green-Schwarz anomaly cancelation mechanism.


\vspace{3mm}
From the target space perspective, these kinds of anomalies
manifest themselves in the fact that the action functional
of the theory, supposed to be a function on configuration space,
happens to be, 
in fact,
 a section of a line bundle. There are
(at least) two reasons why this may happen:
\begin{itemize}
  \item
    the path integral over the fermionic fields is
    to be interpreted not as a function over the configuration
    space of the remaining bosonic fields, but as a section
    in a Pfaffian line bundle over that space
    (reviewed in \cite{FreedIII});
  \item
   the standard action functional for higher abelian gauge
   fields in the presence of electric \emph{and} magnetic 
   charges is also in general just a section of 
   a line bundle over configuration space 
   (discussed in \cite{Freed}).
\end{itemize}
If the tensor product of these two line bundles, namely of the 
Pfaffian and Charge line bundles, 
is a nontrivial
bundle with nontrivial connection then the action is
\emph{anomalous}.
The Green-Schwarz anomaly cancelation mechanism is to 
introduce electric string and magnetic 5-brane charges in 
precisely such a way that the line bundle on configuration
space thus introduced cancels the nontriviality of the
given Pfaffian line bundle due to the fermions in the theory.

\begin{table}[h]
\begin{tabular}{|cc|}
\hline \\
    \begin{minipage}[t]{8cm}
      The complex {\bf anomaly line bundle}
      with connection
      over the space $\mathrm{conf}_{\mathrm{bos}}$ of bosonic fields
      is the tensor product of 
      a Pfaffian line bundle $\mathrm{Pfaff}$ from the
      fermionic path integral
      and another line bundle, $\mathrm{Charge}$,
      due to the presence of 
      electric and magnetic charges.
    \end{minipage}
    &
    \xymatrix{  
       \mathrm{Pfaff}\otimes \mathrm{Charge}  
       \ar[d]
       \\
       \mathrm{conf}_{\mathrm{bos}}
    }
    \\
    \begin{minipage}[t]{8cm}
      The {\bf action funcional} $e^{-S}$ is supposed to be
      a complex function on $\mathrm{conf}_{\mathrm{bos}}$, but
     is in general, in fact, a section of the anomaly line bundle.
    \end{minipage}
    &
    \xymatrix{  
       \mathrm{Pfaff}\otimes \mathrm{Charge}  
       \ar[d]
       \\
       \mathrm{conf}_{\mathrm{bos}}
       \ar@/_1pc/[u]_{e^{-S}}
    }
    \\
    \begin{minipage}[t]{8cm}
      In order for the starting point of quantization to be well
      defined one needs {\bf anomaly cancelation}:
      the anomaly line bundle needs to be
      trivializable and one needs a choice of trivialization
      that identifies it with the trivial line bundle with
      trivial connection.
    \end{minipage}
    &
    \raisebox{-13pt}{
    \xymatrix{  
       \mathrm{Pfaff}\otimes \mathrm{Charge}  
       \ar[r]^{\simeq}
       &
       (\mathrm{conf}_{\mathrm{bos}} \times \mathbb{C},\nabla = 0)
    }
    }
   \\
   \\
    \vspace{4pt}
    \begin{minipage}[t]{8cm}
     The obstruction to this trivialization is the {\bf anomaly}.
     The curvature of the connection on 
     $\mathrm{Pfaff}\otimes\mathrm{Charge}$ is
     called the {\bf local anomaly}, its holonomy the 
     {\bf global anomaly}.
    \end{minipage}
    &
    \raisebox{-13pt}{
    \begin{tabular}{cc}
      local anomaly: & $\mathrm{curv}(\mathrm{Pfaff}\otimes\mathrm{Charge})$
      \\
      global anomaly: &  $\mathrm{Hol}(\mathrm{Pfaff}\otimes\mathrm{Charge})$
    \end{tabular}
    }
    \\
\hline
\end{tabular}           
\vspace{2mm}
  \caption{
    {\bf Anomalies} arising from the fact that the 
    bosonic action functional
    is, a priori, not a function on the bosonic configuration
    space of fields $\mathrm{conf}_{\mathrm{bos}}$,
    but a section of a line bundle over that space. 
  }
\end{table}

\vspace{3mm}
In the following we describe the anomaly cancelation
condition known as the \emph{dual} Green-Schwarz mechanism 
\cite{Gates} \cite{SS} \cite{Freed} and 
related to super 5-branes \cite{LT}. 
It can be obtained
from the worldvolume perspective of the super 5-brane again
as a generalization of how the Spin-condition for the target
of a spinning particle is found. Alternatively, it can be
found from the condition that the index of the total Dirac 
operator (on the fermionic fields called
``dilatino'', ``gravitino'' and ``gaugino'') of the 
effective target space field theory of the heterotic string
vanishes.

\vspace{3mm}
In String theory, the need for String structures was originally found 
in terms of anomaly cancelations, either from the target space
perspective or from the worldsheet perspective:
The effective field theory of the heterotic string on a Spin target space $X$
involves the the (pseudo-)Riemannian metric structure on $X$
with Levi-Civita $\mathrm{SO}(10)$ connection $\omega$ 
(needed for gravity) 
and a gauge bundle $E \to X$ with connection $A$.
    The spinorial field content of the heterotic background theory
    consists of sections of three different bundles: the Spin
    bundle $S$ associated to the principal Spin bundle over 
    the Spin manifold $X$, as well as 
    its tensor products with $T X$ and with the gauge bundle $E$.
    Sections of $S$ are states of the
      the dilatino field, those
    of $S \otimes TX$ correspond to the  
    gravitino field, and those of $S \otimes E$ correspond to the
     the gaugino field. There is a Dirac operator
    associated with all three of these fields, denoted
    $D$, $D_{TX}$ and $D_E$, respectively.
 The anomaly cancellation condition, which ensures that
    the action functional for these fields is a well defined
    function on configuration space, is that a particular linear 
    combination of the
    indices of these Dirac operators vanishes.

\begin{table}[h]
  \begin{center}
  \begin{tabular}{ccccc}
    {\bf Spin bundle}
    &
    {\bf name of field}
    &
    {\bf \begin{tabular}{c}symbol for \\ Dirac operator\end{tabular}}
    &
    {\bf contribution to anomaly}
    \\
    \hline
    &&&
    \\
    $S$ & 
    dilatino
    &
    $D$
    &
    ${\mathcal I}_{1/2}(R) := \mathrm{Index}(D)$& $ = \hat A$
    \\
    \hline
    &&&
    \\
    $S \otimes TX$
    & gravitino&
    $D_{TX}$
    &
    ${\mathcal I}_{3/2}(R) := \mathrm{Index}(D_{TX})$ & $=\frac{1}{8}L(R) 
    + {\mathcal I}_{1/2}(R)$ 
    \\
    \hline
    &&&
    \\
    $S \otimes E$
    & gaugino &
    $D_{E}$
    &
    ${\mathcal I}_{1/2}(R, F) := \mathrm{Index}(D_E)$ & $= \mathrm{ch}(E)\wedge \hat A$
 \\
 \hline
  \end{tabular}
  \end{center}
  \vspace{2mm}
  \caption{
    \label{fermionic fields}
    {\bf Anomaly contributions} 
     from the three different fermionic fields 
     (sections of spin bundles)
     of the target space theory. $L(R)$ is the Hirzebruch $L$-polynomial.
  }
\end{table}

\vspace{3mm}
In the notation of table \ref{fermionic fields},
the anomaly cancelation formula involves  the degree twelve form part of the  identity \cite{AW} \cite{GS} \cite{GSW}
\(
{\mathcal I} := {\mathcal I}_{3/2}(R) - {\mathcal I}_{1/2}(R) 
+ {\mathcal I}_{1/2}(F, R) = 0.
\)
Here $R$ and $F$ are the 
curvature 2-forms of the tangent bundle $TX$ and of the gauge bundle $E$,
respectively.
The second of the three terms in the middle
is ${\rm Index}(D) = \widehat{A}$, the index of the 
uncoupled Dirac operator given in terms of the $A$-genus via the
index theorem. The first term is ${\rm Index}(D_{TM})$, 
 the index of the Dirac operator coupled to the tangent vector bundle, 
i.e. $S \otimes TM$. 
The third term is ${\rm Index}(D_E)$ is the index of the Dirac operator coupled
to the gauge vector bundle, whose curvature is F, i.e. the vector bundle
is ${\rm Spin}(M) \otimes E$.
This is equal to ${\rm ch}(E) \wedge \widehat{A}$, by the index theorem. 

\vspace{3mm}
This fermionic anomaly corresponds to a 
line bundle with connection on configuration space whose 
curvature 2-form is the integral of a 
certain 12-form (see \cite{Sch} for more detail) over 
target space $X$
\(
  \mathrm{curv}(\mathrm{Pfaff}) = - \int_X I_4 \wedge I_8 
  \,.
\) 
A similar integral encodes the curvature of the 
anomaly line bundle due to electric string current $j_E$
(a 4-form) and
magnetic current $j_B$ (an 8-form):
\(
  \mathrm{curv}(\mathrm{Charge}) = \int_X j_E \wedge j_B 
  \,.
\) 
Anomaly cancellation demands that we identify
$I_4$
(right hand of equation (\ref{H3}) below) and 
$I_8$ (right hand of equation (\ref{dH7}) below), respectively, as 
the magnetic currents for the field strengths 
$H_3$ and $H_7$ that appear in the direct and the dual
Green-Schwarz anomaly cancelation conditions.

\subsection{Anomalies in terms of $H_3$ and $H_7$}

Here we summarize the two pictures we have emphasized.

\vspace{3mm}
\noindent {\bf The standard Green-Schwarz mechanism via $H_3$.}
The $H$-fields appear in two theories of interest to us: 
 Type II  
and heterotic theories on ten-dimensional $X$, respectively. 
In type II, the expression involves the Ramond-Ramond (RR) fields.
The direct Green-Schwarz formalism \cite{GS} on a ten-dimensional manifold 
$X$ leads to the appearance of 
the three-dimensional Chern-Simons term via the Chapline-Manton 
coupling \cite{CM} which makes $H_3$ no longer closed.  
Mathematically, this means we assume in that case that
in the heterotic  theory there is the usual $H$-field, a priori an $\R$-valued differential form
--  as presented by supergravity --
which gets modified from being closed to
\(
\frac{1}{2 \pi} dH_3= {\rm ch}_2(A) - \frac{1}{2}p_1 (\omega),
\label{H3}
\) 
where $A$ is the gauge connection for an $E_8 \times E_8$ or ${\rm Spin}(32)/\Z_2$
vector bundle $E$ and $\omega$ is the metric connection, so that $p_1(\omega)$ is 
the first Pontrjagin form of the tangent bundle $TX$.
Recall \cite{Kill} that the right hand side of (\ref{H3})
is the (the image in real cohomology of the) obstruction to 
having a String structure on the 
virtual difference bundle $TM -E$, 
and that (\ref{H3}) therefore says
that this obstruction class has to vanish.
Note the contrast of (\ref{H3}) in heterotic string theory to the case of type II 
string theory on $X$, where $dH_3=0$. 

\vspace{3mm}
\noindent {\bf The dual Green-Schwarz mechanism via $H_7$.}
The $H$-field $H_3$ can be viewed as 
the
dual of a
field $H_7$ where. rationally,
this is just Hodge duality $H_7 := \star H_3$
The `Bianchi identity' of $H_7$ depends on the specific string
theory, i.e. heterotic vs. type II, as was the case for $H_3$. 
The reason we say ``rationally" is because
these fields, like other fields in string theory, can give integral and 
even torsion elements in cohomology in the quantum theory. 
In such a case, appropriate notions
of Hodge duality will be needed to clarify the relationship between 
$H_3$ and $H_7$.  We address this in detail in a separate paper.

\vspace{3mm}
 Let us now also consider type II string theory on $X$ in the absence of 
 any Ramond-Ramond fields. This theory also has
 a degree seven 
 dual, $H_7$, of the $H$-field $H_3$. While  $H_3$ in this case 
 is closed,   $H_7$ is not. Instead, from the dimensional 
 reduction from M-theory, $H_7$ satisfies (cf. \cite{DLM} and
  \cite{FS} \cite{MaS})
\(
  dH_7=\frac{p_2(\omega) - 
\left( \frac{1}{2}p_1(\omega) \right)^2}{48},
\label{H7II}
\)
where $p_i$ are the Pontrjagin classes of the tangent bundle $TX$.
Observe that this only involves the topology of spacetime \emph{without}
any gauge bundles.

\vspace{3mm}
On the other hand, in the heterotic case, we have a principal bundle with connection, 
  the corresponding 
modified Bianchi  identity will fail by the Chern characters 
of $E$. The expression is (see \cite{Freed})
\(
dH_7= 2\pi \left( {\rm ch}_4(A) -\frac{1}{48} p_1(\omega) {\rm ch}_2(A) 
+ \frac{1}{64} p_1(\omega)^2 -\frac{1}{48}p_2(\omega) \right).
\label{dH7}
\)
Here $\omega$ denotes the Levi-Civita connection on
the Spin-lift of the tangent bundle, and $A$ is the
connection on the gauge vector bundle $E$. We interpret this 
``dual Green-Schwarz mechanism'' as saying that $H_7$ 
trivializes the obstruction to having a Fivebrane structure on the
pair ($TX$, $E$)

\vspace{3mm}
\noindent {\bf Observations:}
  It is useful to see how (\ref{dH7}) simplifies in various special cases.

\noindent {\bf 1.} If 
we have a String structure on $T X - E$ coming from String structures on both
bundles separately, in that
$p_1(TX)=0={\rm ch}_2(E)$, 
then, at the level of cohomology, equation (\ref{dH7}) is replaced by 
\(
\left[ {\rm ch}_4(E) -\frac{1}{48} p_2(TX) \right]=0\;.
\)
This holds in general when $X$ is 4-connected, in which case the
cohomology of $X$ is only in degrees 0, 5 and 10. 

\noindent{\bf 2.} 
If in addition to  ${\rm ch}_2(E)=0$ we require that $c_1^2$ and $c_2$ 
be equal to zero, 
then in this case, using 
\(
{\rm ch}_4=\frac{1}{24} \left( c_1^4 -4c_1^2 c_2 + 4c_1 c_3 + 2c_2^2 -4c_4 \right),
\)
we have that ${\rm ch}_4(E)=\frac{-1}{6}c_4(E)$ so that
\( 
 \left[ -\frac{1}{6} c_4(E) - \frac{1}{48} p_2(TX) \right]=0\;.
\label{c4p2}
\)

\noindent {\bf 3.} We can write a formula in terms of either the Chern classes or the Pontrjagin 
classes for both factors in (\ref{c4p2}), thus giving specialization of the 
general formula. 
For this we can consider the 
complexification of the tangent bundle to the ten-dimensional spacetime.
We use 
\(
p_j(TM)=(-1)^j c_{2j}(TM_{\C})
\)
to write the differential form representative of (\ref{c4p2}) as
\( 
dH_7= \frac{-2\pi}{6} \left( c_4(A) + \frac{1}{8} c_4(\omega) \right),
\label{c4c4}
\)
where now it is understood that $c_4(\omega)$ is the fourth Chern 
class of the complexified
 tangent bundle with corresponding connection
$\omega$, for which, with an abuse of notation, we use the same symbol 
as for the real connection. 

\vspace{3mm}
It will follow from our results in section \ref{fivebrane structures} that 
fivebrane structures for a pair $(TM,E)$ with connections $( \omega,A)$ 
   require $\mathrm{ch}_4(A)-\frac{1}{48}p_2(\omega)$
   to vanish on top of the string structure.

\subsection{The Chern-Simons forms}
For a principal $G$-bundle $p:P\to X$, a characteristic class $K_j(F)$, expressed as  a polynomial in the $\frak{g}$-valued curvature 
$F$ of polynomial degree $j$, is closed and 
pulled up to the total space is exact
  \(
K_j(F)=dQ_{2j-1}(A,F),
\)
where $Q_{2j-1}(A,F) \in \frak{g} \otimes \Lambda^{2j-1}(M)$ is a ` Chern-Simons 
form'  for $K_j(F)$. This applies to the Chern character as well as the Pontrjagin 
classes. Thus, we can  use this to solve equation (\ref{dH7}) 
on the total space. We denote a specific choice (e.g. by the specific homotopy in the remark below)
as $CS_7(A, F) .$
Using the 
 expressions
\bea
{\rm ch}_4(F)=dCS_7(A, F) 
\\
{\rm ch}_4(R)=dCS_7(\omega, R),
\eea
equation (\ref{c4c4}) becomes 
\(
dH_7= 2 \pi \left[ dCS_7(A, F) + \frac{1}{8} dCS_7(\omega, R)  \right].
\)
This means that 
 $H_7$ can be taken to be of the form 
\(
H_7= 2\pi \left( CS_7(A, F) + \frac{1}{8} CS_7(\omega, R) \right). 
\)
We view this setting of $H_7$ equal to the degree seven Chern-Simons forms 
as the degree seven analog of the Chapline-Manton coupling.

\vspace{3mm}
\noindent {\bf
Remarks.}
{\bf 1.} A specific formula for the Chern-Simons form corresponding to the Chern character
 can be obtained 
  by using the
   homotopy formula \cite{CS} (or in the physics
 literature \cite{ZWZ} \cite{AG})
 \(
 CS_{2j-1}(A, F) = \frac{1}{(j-1)!}\left( \frac{i}{2\pi} \right)^j
 \int_0^1 dt ~ {\rm Str} (A, F_t^{j-1}),
 \)
 where Str is the symmetrized trace and $A_t$ is a connection that interpolates 
 between connections $A_0$ 
 and $A_1$
 \(
 A_t=A_0 + t(A_1 - A_0),
 \)
  with corresponding curvature 
 \begin{eqnarray}
 F_t &=& dA_t + A_t^2
 \nonumber\\
 &=&t dA + t^2 A^2
\nonumber\\
&=& tF + (t^2-t)A^2.
 \end{eqnarray}
 For $j=4$, 
\(
CS_7(A, F) = \frac{1}{6 (2\pi)^4} \int_0^1 dt~ {\rm Str} \left(A, F_t^3  \right),
\)
where $F_t$ is the curvature of the connection $A_t=tA$ that 
interpolates between the zero connection at $t=0$ and the connection $A$ at $t=1$. 
An analogous formula holds for 
$CS_7(\omega, R)$. 

\noindent {\bf 2.} Unlike the Chern character, the Chern-Simons form is not gauge-invariant. 
Under a transformation $A \to A^g= g^{-1} (A + d) g$, 
\(
CS_7(A^g, F^g) - CS_7(A, F) = -\frac{3!}{7!} {\left( \frac{i}{2\pi} \right)}^4 {\rm tr} 
\left[ (g^{-1}dg)^7 \right] + d \beta_6,
\)
where $\beta_6$ is a six-form which can be chosen 
by applying the chosen homotopy operator on
the gauge-transformed Chern-Simons form. 
 
\noindent {\bf 3.}
Alternatively, we see that (\ref{c4p2}) is obtained from (\ref{dH7}) by
setting all \emph{decomposable} characteristic forms 
(all nontrivial wedge products of two characteristic forms) to 0, 
which is
the same as saying that all characteristic classes \emph{suspending to 0}
are set to 0.
Recall that a characteristic form $P(F_A)$ on a $G$-bundle $p : P \to X$ is
said to suspend to the form $\mu(i^* A)$ on $G \stackrel{i}{\hookrightarrow} P $ 
if there is a form 
$\mathrm{CS}_P(A,F_A)$ on $P$ such that
\(
  d \mathrm{CS}_P(A,F_A) = p^* P(F_A)
\)
and
\(
  i^* \mathrm{CS}_P(A,F_A) = \mu(i^* A)
  \,.
\)
Put more simply, recall that a  form $\omega\in H^*(X)$ on the base of a bundle 
$p : E \to X$ is said to suspend to the form $\mu$ 
on $F \stackrel{i}{\hookrightarrow} E $ 
if there is a form $\nu$ on $E$ such that
$d \nu = p^*\omega$
and
$
  i^* \nu = \mu
$.
A decomposabe characteristic form $P(F_A) \wedge P'(F_A)$ 
necessarily suspends to a 0-form, since
\(
  p^*(P(F_A) \wedge P'(F_A)) = 
  d \mathrm{CS}_{P\wedge P'}(A, F_A))
  :=
  d
  \left(
    \mathrm{CS}_P(A,F_A) \wedge P'(F_A)
  \right)
\)
and since 
\(
  i^*
  \left(
    \mathrm{CS}_P(A,F_A) \wedge P(F_A)
  \right)
  = 0,
\)
because $i^* F_A = 0$, and hence $i^* P'(F_A) = 0$.

\section{Fivebrane Structures}

\label{fivebrane structures}

We recall how String structures are lifts of Spin bundles through the 
3-connected cover of $\mathrm{Spin}(n)$ and how this lift is obstructed by the
fractional first Pontrjagin class called $\frac{1}{2}p_1$. 
Then we define \emph{Fivebrane structures}
as lifts of the resulting String bundles through the 7-connected cover of
$\mathrm{Spin}(n)$.
We demonstrate that this lift is obstructed by a degree 8 characteristic class
of String bundles, which is a fractional second Pontrjagin class, $\frac{1}{6}p_2$,
of the String bundle.

\subsection{The homotopy groups of $SO(n)$}

\label{The homotopy groups of SO}
 
The homotopy groups of the orthogonal group $O(n)$, for $n$ sufficiently large, 
are 
\(
\pi_k \left( O(n) \right)= \left\{ \begin{array}{ll}
\Z_2 & {\rm for}~k=0,1~{\rm mod}~8 \\
\Z & {\rm for}~k=3,7~{\rm mod}~8 \\
0 & {\rm otherwise}
\end{array}
\right.
\label{pi}
\,.
\)
See also the first two rows of figure \ref{homotopy groups of O(n)}.
The condition on $n$ is best understood by considering the
\emph{stable orthogonal group}, also know as the infinite orthogonal group, 
which is defined as the 
direct limit of the sequence of inclusions
\(
O(1) \subset O(2) \subset \cdots \subset O=\bigcup_{k=0}^{\infty} O(k).
\)
That the homotopy groups of $\mathrm{O}(n)$ stabilize, i.e. that 
for $n>k+1$, one has $\pi_k(\mathrm{O})=\pi_k(\mathrm{O}(n))$,
follows from the fact that 
the inclusion $O(n) \hookrightarrow O(n+1)$ 
is $(n-1)$-connected (see below for notation) as we have the fiber bundle
\(
O(n) \to O(n+1) \to S^n,
\)
i.e. $S^n$ is the homogeneous space $O(n+1)/O(n)$. 

\vspace{3mm} 
Below the stable range, i.e. for $n \leq k +1$, the description 
of the homotopy groups of $\mathrm{SO}(n)$ becomes 
incomplete because 
one is essentially looking at the homotopy groups of spheres, which are
not completely known. For example,
the homotopy groups of $\mathrm{SO}(3)$ are the same as the homotopy groups of
the sphere $\mathrm{S}^3$ which are known only in specific degrees (but nonetheless
in a range sufficient for most practical purposes).
But in the applications of all these considerations to string theory which we have
in mind,  the base manifold $X$ from which one wishes to 
consider homotopy classes of maps into $B (\mathrm{SO}(n)\langle k \rangle)$
is $(n=10)$-dimensional. Therefore in this application figure 
\ref{homotopy groups of O(n)} is fully applicable.
One obtains a sequence of topological groups from $\mathrm{O}(n)$ by successively
passing to its $k$-connected covers. This process is known as the 
``killing of homotopy groups''. See the third row of figure 
\ref{homotopy groups of O(n)}.

\vspace{3mm}
The description of the unitary group is analogous but much simpler due to the fact 
that Bott periodicity in this case has period two instead of period eight. In the stable 
range, i.e. for $i<2n$, the inclusion
$U(n)$ into $U(n+1)$ induces an isomorphism between the homotopy
groups $\pi_i(U(n))$ and $\pi_i(U(n+1))$. The infinite unitary group $U$ is defined
in an analogous way as the orthogonal group above. 
The Bott periodicity theorem implies that the homotopy groups of 
$U$ are particularly simple: 
$\pi_i(U)$ is trivial if $i$ is even and isomorphic to $\Z$ if $i$ is odd.

\vspace{3mm}
The standard notation 
for the $k$-connected cover of a space $X$
for which all the homotopy groups vanish up to \emph{and including} $\pi_k$
is  $X\langle k+1 \rangle$.
Thus the groups $U\langle 2k \rangle$ and $O\langle 2k \rangle$ denote 
the connected covers of the unitary and orthogonal group, respectively,
having the first potentially nontrivial homotopy group in dimension $2k$.
	For example, $O\langle 3\rangle$
refers to the orthogonal group having first homotopy group in dimension three, which 
means that the first three homotopy groups are killed (starting with $\pi_0$) and therefore is 
2-connected. 
The result of killing the first two homotopy groups of $O(n)$
are very familiar: these are the groups $SO(n)$ and $\mathrm{Spin}(n)$,
respectively.
The  group $SO(n)$ is the connected component of the  group $O(n)$  
and $\mathrm{Spin}(n)$
is the simply-connected cover of $SO(n)$.  
Less familiar but by now well known is the result of killing the next nontrivial
homotopy group beyond that, $\pi_3$: this is the String group $\mathrm{String}(n)$.

\vspace{3mm}
It is noteworthy that $O(n)$, $SO(n)$ and $\mathrm{Spin}(n)$
are not just topological groups, but of course carry the structure of a Lie group.
On the other hand, there is no known model for $\mathrm{String}(n)$
as a group with a manifold structure. Notice
that since all \emph{compact} Lie groups $G$ have $\pi_3(G) = \mathbb{Z}$ or several copies of 
$\mathbb{Z}$, 
$\mathrm{String}(n)$ has no chance of being a compact Lie group. 
The known models for $\mathrm{String}(n)$ are ``huge'' topological spaces. 
Notice however that \cite{BCSS} give a Lie model 
for $\mathrm{String}(n)$ when
regarded not as a mere group, but as a Lie 2-group 
(a 2-group whose space of objects and of morphisms is
a Fr{\'e}chet) manifold.

\vspace{3mm}
\noindent {\bf
Caveat: notation for higher connected covers.}
A comment on the notation used for the connected covers of the corresponding classifying spaces
is in order.
 In another
approach,  the connected covers of the groups are defined not directly 
but via the corresponding classifying space 
\(
G \langle n \rangle :=\Omega B G \langle n \rangle,
\label{bad}
\)
where
$\Omega$ denotes the based loop space.
However, for our purposes we are also interested in the classifying space itself. To avoid confusion
(hopefully),
we will use the notation  
\(
B( G\langle n \rangle) = (BG)\langle n+1 \rangle.
\)
Note the shift in $n$ compared to (\ref{bad}).
 For example, $(BO)\langle 4\rangle = B(O\langle 3\rangle)$
refers to the classifying space having first homotopy group in dimension four, which 
means that the first four  homotopy groups are killed, and therefore is 
3-connected.

\subsection{String structures revisited}
\label{String structures revisited}

Recall how the group $\mathrm{String}(n)$ can be constructed 
from $\mathrm{Spin}(n)$ via the sequence 
\(
1 \to K(\Z, 2) \to  \mathrm{String}(n) \to \mathrm{Spin}(n)
\to 1,
\label{kz2}
\)
where ${\mathrm{Spin}}(n)$ is the simply connected cover of ${\mathrm SO}(n)$ 
and, in fact, the first nonzero homotopy group of the Lie 
group ${\mathrm Spin}(n)$ is $\pi_3 = \Z$.
It follows from the Hurewicz theorem that the first nonzero cohomology group occurs
at the same degree, i.e. $H^3(\mathrm{Spin}(n),\Z) = \Z$. 
This singles out the (homotopy class of a) map
\(
  f : \mathrm{Spin}(n) \to K(\Z,3) \simeq B K(\Z,2)
\label{bkz2}
\)
that generates $H^3(\mathrm{Spin}(n),\Z)$ under the identification
$H^3(X, \Z) \simeq \left[X, K(\Z, 3) \right] $, the set of homotopy classes of maps $X\to K(\Z, 3) $. 
This map
classifies a principal $K(\Z,2)$-bundle over $\mathrm{Spin}(n)$
\(
 \raisebox{30pt}{
 \xymatrix{
    **[r] \hspace{-15pt}\mathrm{String}(n) := f^* E K(\Z,2) \ar[rr] 
    \ar[dd]
     && 
     E K(\Z,2)
     \ar[dd]
    \\
    \\
    \mathrm{Spin}(n) \ar[rr]^<<<<<<<<<f &&  **[l] K(\Z,3)\simeq B K(\Z,2)\hspace{-15pt}
  }
  }
\)
and this bundle is the extension (\ref{kz2}).
Applying the classifying functor on (\ref{kz2})  leads to a  weakly exact (homotopy exact) sequence
\(
K(\Z, 3) \to B{\rm String}(n) \to B{\rm Spin}(n).
\)
The corresponding complex analog is 
\(
K(\Z, 3) \to (BU)\langle 6 \rangle \to BSU
\)
and on into $K(\Z, 4)$ obtained by mapping from $BSU$ 
by the Chern class $c_2$. The structures are summarized in the following diagram
\(
  \raisebox{60pt}{
  \xymatrix{
    (BU)\langle6\rangle
    \ar[rr]
    \ar[dd]
    &&
    **[l]B \mathrm{String}= (BO)\langle 8\rangle \hspace{-30pt}
    \ar[dd]
    \\
    \\
    BSU=(BU)\langle 4 \rangle 
    \ar[rr]
    \ar[dd]
    &&
    B \mathrm{Spin} = (BO) \langle 4 \rangle
    \ar[dd]
    \ar[rr]^{\frac{1}{2} p_1}
    &&
    K(\Z,4)
    \\
    \\
    BU 
    \ar[rr]
    \ar[rr]
    \ar[dr]
    &&
    B SO = (BO) \langle 2 \rangle
     \ar[dd]
    \ar[rr]^{w_2}
    &&
    K(\Z_2, 2)
    \\ 
    &
    K(\mathbb{Z},2)
    \ar[urrr]|<<<<<<<<<<<<<<<<<<<<<<{\makebox(12,12){}}_<<<<<<<{\mathrm{mod}\, 2}
    \\
    &&
    BO=(BO) \langle 1\rangle
  \ar[rr]^{w_1} 
  &&
  K(\Z_2, 1).
  }
  }
\label{web1}
\)

\vspace{3mm}
\noindent {\bf
Remarks.}
{\bf 1.} $(BU)\langle6\rangle$ and $(BO)\langle 8\rangle$
are the homotopy fibers of the respective maps to $K(\Z,4)$.

\noindent {\bf 2.} The map from $B{\rm SU}$ to $B{\rm Spin}$ is an isomorphism 
on fourth cohomology $H^4$. This comes from the isomorphism 
\(
{\rm Spin}(3) \cong SU(2) \cong S^3.
\) 

\noindent {\bf 3.} The composite map from $B{\rm SU}$ to $K(\Z, 4)$ in the 
second row is given by the second Chern class $c_2$.

\noindent {\bf 4.} The composite map from $B{\rm U}$ to $K(\Z_2, 2)$ in the
third row is given by $c_1$ followed by reduction mod 2. 
That is, we have the following commutative
diagram
\(
  \xymatrix{
BU 
\ar[rr]^{c_1}
\ar[dd]
 &&
 K(\Z, 2)
 \ar[dd]^{{\rm mod}~2}
\\
\\
BU
\ar[rr]^{w_2}
&&
K(\Z_2, 2).
}
\)

\subsection{Cohomology of the connected covers}

\subsubsection{In characteristic 0}
In this section we work in characteristic zero and we (briefly) address
torsion in (\ref{torsion}). 
The definition of $\mathrm{String}(n)$ as  the extension of $\mathrm{Spin}(n)$ 
in (\ref{kz2}) induced by the map (\ref{bkz2})
allows us to compute the cohomology of $\mathrm{String}(n)$. It is easier 
to compute first the cohomology of $B{\rm String}(n)$ from the 
homotopy fibration sequence
\(
K(\Z, 3) \to B{\rm String}(n) \to B{\rm Spin}(n),
\)
since in characteristic 0 we have (all coefficients are rational)
\(
H^*(K(\Z, 3) ) \simeq H^*(S^3)
\) 
and we can apply the  long exact Gysin sequence
\(
\cdots \longrightarrow H^n(E)\longrightarrow^{\!\!\!\!\!\!\!\!\!\!\pi_*}~~H^{n-k}(M)\longrightarrow^{\!\!\!\!\!\!\!\!\!\!e\wedge}~~H^{n+1}(M)\longrightarrow^{\!\!\!\!\!\!\!\!\!\!\pi^*}~~H^{n+1}(E)\longrightarrow \cdots
\)
where $e\wedge$ is the wedge product of a differential form with the Euler class $e.$
The result is 
\(
H^*(B\mathrm{String})\simeq P[p_2,p_3,\dots],
\)
from which it readily follows that $H^*(\mathrm{String})$ is an exterior algebra on generators
$x_i$ of degree $4i-1$ starting with $i=2.$

\vspace{3mm}
A very close analog of the above can be carried over starting with $BU.$ The result is 
\(
H^*((BU) \langle 2k \rangle) \simeq P\left[ c_k, c_{k+1}, \cdots \right]/ (c_k) = P\left[ c_{k+1}, \cdots \right].
\label{rationalBU}
\) 
 Except for the change in indexing, the proof is the same as 
 for $B\mathrm{String} = B\mathrm{Spin} \langle 4\rangle.$

\vspace{3mm}
\noindent {\bf
Remarks.}
{\bf 1.} The rational 
cohomology of $(BU) \langle k \rangle$ can be calculated 
in several ways, for example using the Gysin sequence, as we saw 
above, or 
by induction on $k$ using the Leray-Serre spectral sequence. 
While the first method is more straightforward, we include the
second one in  Appendix B because it is likely to be needed for applications 
in future work  where  torsion is  important.

\noindent {\bf 2.} In the geometric description of $K(\Z, 2)$ as $PU(\mathcal{H})$, the projective
unitary group on an infinite-dimensional Hilbert space $\mathcal{H}$, 
the 
$H$-field $H_3$ occurs as the canonical degree three
class $x_1$ of $PU(\mathcal{H})$
 bundles. In the operator algebra language this is called the Dixmier-Douady class
 \cite{DD}.
  From a physical point of view if $X$ is ten-dimensional, it is  natural to ask whether a
 similar interpretation
of the dual field $H_7$ can be given. From a mathematical point of view, it is natural to ask 
whether a similar construction to the degree three case 
can be carried out with the first non-trivial 
homotopy group $\pi_7$ of $\mathrm{String}$.

 \subsubsection{Over the integers and the integers mod $p$.}
\label{torsion}


{\bf 1.} The cohomology for all connected covers for the classifying spaces 
of the orthogonal and the 
unitary groups are also known.
At $p=2$ this is a result of Stong \cite{Stong}. For general 
$p$, this is a result of 
W. Singer \cite{Singer}, in terms of generators and relations, 
and Giambalvo \cite{Giam}.

\noindent {\bf 2.} $H_{*}\left( (BU)\langle 6 \rangle; \Z  \right)$ is calculated as a ring 
in \cite{AHS}. 
  It is torsion-free and concentrated in even degrees.
  From this 
the cohomology can be read via 
\(
H^*= {\rm hom}\left( H_*((BU)\langle 6 \rangle; \Z), -  \right),
\)
as a map from rings to sets. This is related in an interesting way
to the what are known as  cubical structures on the additive group \cite{AHS}.

\subsection{Fractional Pontrjagin classes}

When our manifold $X$ has extra structure, such that it admits a lift 
of the structure group of its 
tangent bundle to a higher connected cover of $O(n)$, we can refine the
Pontrjagin classes of $X$ to \emph{fractional} Pontrjagin classes.
For $X$ a $d$-dimensional orientable manifold, its $k$th Pontrjagin class 
$p_{k}$ is taken to be the Pontrjagin class of the tangent bundle $TX$, 
regarded as an associated $SO(d)$-bundle.
So it is given by a map
\(
  \xymatrix{
    X \ar[rr]^f
    &&
    B SO(d)
    \ar[rr]^{p_k}
    &&
    K(\mathbb{Z},4k)
  }
  \,.
\)

\subsubsection{Spin structures and the first Pontrjagin class}
\label{spin structures and first pontrjagin}

Saying that a $d$-dimensional manifold $X$ is \emph{spin} means we can lift the classifying map
$f$ of its tangent bundle $TX$ to a map
$\hat f : X \to B \mathrm{Spin}(d)$. Since $\mathrm{Spin}$
is 2-connected and  $\pi_3(\mathrm{Spin}) = \mathbb{Z}$,
it follows that $H^4(B\mathrm{Spin},\mathbb{Z}) = \mathbb{Z}$.
If we denote the generator of this fourth cohomology group by
$\omega,$ the situation looks as follows:
\(
  \raisebox{70pt}{
  \xymatrix{
    &&
    B \mathrm{Spin}(d)
    \ar[dd]^\pi
    \ar[rr]^{\omega}
    &&
    K(\mathbb{Z},4)
    \\
    \\
    X \ar[rr]^f
    \ar[uurr]^{\hat f}
    &&
    B \mathrm{SO}(d)
    \ar[rr]^{p_1}
    &&
    K(\mathbb{Z},4)
  }
  }
  \,.
\)
Since, by assumption, $\omega$ is a generator, it must be true that
$\pi^* p_1$ is an integral multiple of $\omega$. 
One finds that  $\omega$ can be chosen so that
\(
 \pi^* p_1 = 2 \omega
 \,,
\)
i.e.
\(
  \raisebox{70pt}{
  \xymatrix{
    &&
    B \mathrm{Spin}(d)
    \ar[dd]^\pi
    \ar[rr]^{\omega}
    &&
    K(\mathbb{Z},4)
    \ar[dd]^{\cdot 2}
    \\
    \\
    X \ar[rr]^f
    \ar[uurr]^{\hat f}
    &&
    B \mathrm{SO}(d)
    \ar[rr]^{p_1}
    &&
    K(\mathbb{Z},4)
  }
  }
  \,.
\)
This motivates the notation
\(
  \omega := \frac{1}{2}p_1
\)
for the generator of $H^4(B \mathrm{Spin},\mathbb{Z})$.
Accordingly, the pullback
\(
  \frac{1}{2}p_1(X) := {\hat f}^* \frac{1}{2} p_1 :
  \xymatrix{
    X \ar[r]^{\hat f} & B \mathrm{Spin} \ar[r]^{\frac{1}{2}p_1} & K(\mathbb{Z},4)
  }
\)
is ``half the Pontrjagin class''  of  the spin-manifold $X$.
Notice that for $\frac{1}{2}p_1(X)$ to be zero, the vanishing of $p_1(X)$ is a necessary but 
not a sufficient condition: $\frac{1}{2} p_1(X)$ might be non-vanishing but  2-torsion.

\subsubsection{String structures and the second Pontrjagin class}
\label{String structures and the second Pontrjagin class}

The same kind of reasoning continues to apply as we keep killing 
homotopy groups of $O(d)$. 
We say that $X$ 
admits a \emph{string} structure or that $X$ is string 
if the classifying map $f$ for $TX$
lifts to $B \mathrm{String}(d).$
Now $B\mathrm{String}$ is 7-connected and  
$H^8(B\mathrm{String},\mathbb{Z}) \simeq \mathbb{Z}$. Let $\nu$
denote the corresponding generator. The pullback 
$\pi^* p_2$ of the second
Pontrjagin class
 has to be an integer multiple of this generator. 
In the next section (see Proposition \ref{p26} section \ref{5b}), we will 
show that  the integer multiple is 6:
\(
  \raisebox{70pt}{
  \xymatrix{
    &&
    B \mathrm{String}(d)
    \ar[dd]^\pi
    \ar[rr]^{\nu}
    &&
    K(\mathbb{Z},8)
    \ar[dd]^{\cdot 6}
    \\
    \\
    X \ar[rr]^f
    \ar[uurr]^{\hat f}
    &&
    B \mathrm{SO}(d)
    \ar[rr]^{p_2}
    &&
    K(\mathbb{Z},8)
  }
  }
  \,.
\)
Therefore we should give the generator $\nu$ the name
$p_2/6$
and define
\(
  \frac{1}{6}p_2(X) := \hat{ f}^* \frac{1}{6} p_2,
\)
the fractional second Pontrjagin class of the String manifold $X$.
 Later (see equation (\ref{qs})) we will make
use of the spin characteristic classes, which better describe 
spin bundles than do Pontrjagin classes.

\subsection{Fivebrane Structures}
\label{5b}

The point of view we take is that in the same way that $H_3$ was part of 
the obstruction to lifting a Spin bundle to a String bundle, 
 we would like to interpret $H_7$ as another obstruction.
Further, $H_7$ serves as a higher obstruction in the sense that it 
makes sense to talk about it once the `lower' obstructions vanish. 
The next task is to make this more precise. 

\vspace{3mm}
The first nonzero homotopy group of the topological
group $\mathrm{String}(n)$ is
$\pi_7\simeq \mathbb{Z}$.
Then, again, the Hurewicz
theorem implies that
the first nonzero cohomology group occurs in degree 7. 
As 
$H^7(X, \Z)=\left[ X, K(\Z, 7) \right]$ and $K(\Z, 7)= B K(\Z, 6)$,
it follows that $K(\Z, 6)$ is a fiber in a 
nontrivial fibration  
sequence with $\mathrm{String}(n)$ 
as the base. From the structure of the homotopy groups (\ref{pi}),
the extension will be $O\langle 8 \rangle$. Thus (compare also Def. \ref{maindef}) 

\begin{definition}
The extension
$
1 \longrightarrow K(\Z, 6) \longrightarrow O\langle 8 \rangle \longrightarrow  
  \mathrm{String} 
\longrightarrow 1
$,
classified by the canonical map $\mathrm{String} \to K(\mathbb{Z},7)$,
which replaces the sequence (\ref{kz2}), we call the 
$\mathrm{Fivebrane}$-extension
$
  \mathrm{Fivebrane} := O\langle 8 \rangle
  $.
  \label{def 2}
 \end{definition}

\vspace{3mm}
\noindent {\bf Remarks.} 
{\bf 1.} The fact that $\mathrm{String}$ occurs as the base of the sequence in Def. \ref{def 2}
is compatible with the interpretation of the classes of $H_3$ and $H_7$ as 
corresponding to the fibrations being pulled back from $K(\Z, 3)$ and $K(\Z, 7)$ respectively. 

\noindent {\bf 2.} When the space is $\mathrm{Spin}$, 
the first Pontrjagin class $p_1$ is divisible by $2$, and the 
obstruction to lifting $\mathrm{Spin}$ to $\mathrm{String}$ is 
$\frac{1}{2}p_1$. When the space is String, the second
Pontrjagin class $p_2$ is divisible by $6$, and the obstruction to 
lifting $String$ to 
$\mathrm{Fivebrane}$ is $\frac{1}{6}p_2$.

\begin{proposition}
The obstruction to lifting a $String$ bundle to an $O \langle 8 \rangle$ bundle is given
by $\frac{1}{6}p_2$. 
\label{p26}
\end{proposition}
 
\proof
 The classifying spaces of the above sequences  induce a map from 
$B\mathrm{String} \to B\mathrm{Spin}$. The composite map 
$B\mathrm{String} \to B\mathrm{Spin} \to BSO$
will map $p_2$ to a multiple of the generator of 
$H^8(B\mathrm{String},\mathbb{Z})$. 
This multiple is obtained 
by noticing that the map from $B\mathrm{String}$ to 
$(BU)\langle 7 \rangle$ 
is an isomophism on $\pi_8$ and that the fourth Chern class 
$c_4$ restricts to 6 times the generator of $\pi_8$ on the complex side. 

\vspace{3mm}
The complex vector bundle on $S^{2k}$ that represents the 
generator of $\pi_{2k}(BU)$
is the $k$-fold (external) tensor product of $(1-L)$, where $1$ is 
the tautological (trivial)
line bundle and $L$ is the Hopf line bundle
on $S^2$. All that 
needs to be done is find the value of $c_k$ on this bundle. Using the
multiplicative properties of the Chern character, this is 
just $(k-1)!$ (see e.g. \cite{Peterson}, Theorem 5.1). Thus for 
$k=4$, the answer is 6 as claimed. 
\endofproof

\vspace{3mm}
It follows from connectivity that there is a commutative pullback diagram 
$(BO)\langle 9 \rangle \to (BU)\langle 10\rangle$ 
sitting over $(BO)\langle 8\rangle \to (BU)\langle 8\rangle$. The map 
$(BU)\langle 10 \rangle \to (BU)\langle 8 \rangle$ is the fibration corresponding to 
$\frac{1}{6}c_4$, which 
restricts to $\frac{1}{6}p_2$ in $(BO) \langle 8\rangle$. The map 
$(BO)\langle 8 \rangle \to (BU)\langle 8\rangle$ is an isomorphism on $\pi_8$ and 
lifts to $(BO) \langle 9 \rangle \to (BU)\langle 10 \rangle$, which is also an isomorphism
on $\pi_8$. 
Thus we have the following diagram, where $B{\rm Fivebrane}$ is as in Def. 1,
\(
  \raisebox{130pt}{
  \xymatrix{
   B{\rm Fivebrane}=(BO)\langle 9 \rangle
    \ar[rr]
    \ar[dd]
    &&
    **[l] (BU)\langle10\rangle \hspace{-10pt}
    \ar[dd]
    \\
    \\
    B{\rm String}=(BO)\langle 8 \rangle 
    \ar[rr]
    \ar[dd]
    &&
    (BU)\langle 8\rangle
    \ar[dd]
    \ar[rr]^{\frac{1}{6}c_4}
    &&
    K(\Z,8)
    \\
    \\
    B{\rm Spin}=(BO) \langle 4 \rangle
    \ar[rr]
    \ar[dd]
    &&
    B SU = (BU) \langle 4 \rangle
     \ar[dd]
      \ar[rr]^{c_2}
    &&
    K(\Z,4)
    \\ 
    \\
   BSO=(BO)\langle2\rangle
   \ar[rr]
   \ar[dd]
   &&
   BU
   \\
   \\
    BO=(BO) \langle 1\rangle
  }
  }
  \label{b5}
\)

\vspace{3mm}
\noindent {\bf
Remarks.}
{\bf 1.} In order for a bundle to be lifted from 
$\mathrm{String}$ to $O \langle 8 \rangle,$ it has  
already to be $\mathrm{Spin}$. 
Thus the first Pontrjagin classes in expressions (\ref{H7II}) 
and (\ref{dH7})
are set to zero and then 
$dH_7$ in both cases is some linear combination of the second Pontrjagin classes,
corresponding to the (lifts of the) tangent and the gauge bundles, $TX$ and $E$
(cf. \cite{Kill}).

\noindent {\bf 2.} That the first column in (\ref{b5}) is real and the second
column is complex follows the complexification map $KO \to KU$
on bundles, since we are interested in the divisibility properties of the
Pontrjagin classes, and those are considered as the pullback of
the complexification, i.e. the Chern classes.
  
\noindent {\bf 3.} The map from $BO$ to $BU$ is an isomorphism on $\pi_{8k}$ and is multiplication
by 2 on $\pi_{8k+4}$. In particular, the map is an isomorphism on $\pi_8$ and 
is multiplication by 2 on $\pi_4$ (see e.g. \cite{Rud}).

\noindent {\bf 4.} In the second row, we could have written $(BU)\langle 6 \rangle$ in the second
entry. The reason for having $(BU)\langle 8 \rangle$ instead is because we would 
then have the map from the real to the complex side, $(BO) \langle 8 \rangle
\to (BU)\langle 8 \rangle$ to be an isomorphism on $\pi_8$. In contrast, the
map $B{\rm Spin} \to BSU$ is not an isomorphism on $\pi_4$ but is in fact
given by multiplication by 2 (by part 3 above) -- see \cite{Rud} for instance. 

\noindent {\bf 5.} In the first row,
we wrote $(BU)\langle 10 \rangle$ instead of $(BU)\langle 9\rangle$ because
there are no homotopy groups of odd degree in $BU$, i.e. killing the homotopy
in degree eight automatically gets us to degree ten.

\noindent {\bf 6.} In both cases, $B{\rm String}$ and $B{\rm Fivebrane}$, we are killing a $\Z$ in 
the homotopy groups.

\noindent {\bf 7.} The map from $(BO)\langle 8 \rangle$ to $(BU)\langle 8 \rangle$ in the second
row is an isomorphism in degree eight because the generator of $\pi_8$
for both spaces is $v_1^4$, where $v_1$ is the Bott generator whose degree
is two. 

\noindent {\bf 8.} The map from $B{\rm Spin}$ to $BSU$ in the third row is given by 
multiplication by two.

\vspace{3mm}
\noindent {\bf
Real vs. Complex.}
In diagram (\ref{web1}) we had the complex spaces in the first 
 column and the real spaces in the second column. In diagram
 (\ref{b5}) we had instead the real spaces in the first column and 
 the complex spaces in the second column. We describe this
 further here. Consider the composite map  
 from $BSU$ to itself factoring through $B{\rm Spin}$
 \(
 \xymatrix{
 BSU \ar[rr]^{\cong~{\rm in~deg}~4} 
 &&
 B{\rm Spin} 
 \ar[rr]^{\times 2}
 &&
 BSU. }
 \)
 The composition is given by multiplication by 2 in degree 4 and
 acts on vector bundles via
 \(
 V \mapsto V \oplus \overline{V}=2V,
 \)
 so that $c_2(2V)=2c_2(V)$. Note that since we have $SU$ bundles 
 then $c_1(V)=0$. 

\vspace{3mm}
 Next going from $BU$ to $BSO$ we have
 \(
 \xymatrix{
 \Z \times BU \ar[rr]^{ \times 2, ~\cong~{\rm in ~deg}~0}_{\rm forget} 
 &&
\Z \times BSO
\ar[rr]^{\cong~{\rm deg}~0}_{\rm complexify} 
&&
\Z \times BU 
 \ar[rr]
 &&
 \Z \times BSO. }
 \)
 The first map amounts to forgetting the complex structure and the 
 second to complexifying. 
 The map from second to the fourth term 
 is $ V \mapsto 2V$, and that from the first to the third is 
 $V \mapsto V \oplus \overline{V}$.

 \subsection{Congruence}
\label{congruence}

\vspace{3mm}
\noindent {\bf
Remarks.}

{\bf 1.} At the integral level,
 there is a very crucial difference between $p_2$ being zero
and $\frac{1}{6}p_2$ being zero, due to the possible existence of very important 
2- and 3-torsion. In other words, $\frac{1}{6}p_2=0$ certainly implies that $p_2$ is
zero, but the converse is not necessarily true -- and in fact in most interesting
cases it is not true. This gives us the important conclusion that: {\it Unlike the 
$String$ case, both two-  and three-torsion are important for $(BO) \langle 9\rangle$
structures.}

\noindent  {\bf 2.} Note that  
$\pi_* O$ has only 2-torsion so that in the stable
range there is no $4$-torsion, so  there is no difference between 
 $\frac{1}{2}p_1$ and $\frac{1}{4}p_1$ as obstructions.
 The latter
is the shift in the quantization condition of the field strength 
$G_4$ in M-theory \cite{Flux} (cf. footnote 2).

\noindent {\bf 3.} Likewise, in the fivebrane case -- again assuming the stable range-- 
there is no difference, as far as the
type of torsion is concerned, between the obstructions $\frac{1}{6}p_2$
and $\frac{1}{48}p_2$. This is because they are both of the form $\frac{1}{2^n3^m}p_2$,
i.e. involve $2-$ and 3-torsion. This leads to the following conclusion:
{\it In the stable range (hence no 8-torsion) for the tangent bundle, our definition 
of the fivebrane structure captures 
the physical condition for a fivebrane structure, i.e. that $\frac{1}{6}p_2$
is essentially the same as $\frac{1}{48}p_2$ in the sense that has just 
been explained.} 

\noindent {\bf 4.} In both case we can avoid subtleties of division by 2 and 3, by {\it inverting}
those two primes. Thus if we invert 2 in the string case, starting with a ring 
$R$ (e.g. $\Z$) we can work with 
the ring $R[\frac{1}{2}]$, and if we invert $6=2\times 3$ in the fivebrane case, we can
work with the ring $R[\frac{1}{6}]$.

\vspace{3mm}
\noindent {\bf
The division by 8.}
We have seen that the obstruction to the fivebrane structure is given
by $\frac{1}{6}p_2$, while the expression appearing in the anomaly is 
$\frac{1}{48}p_2$.
Here we give a description of this further division by 8, which, by the above
remark, does not generate any new torsion. In particular this means that
the situation in the fivebrane case is better than that in the string case, where
 there was a crucial difference between $p_1$ vanishing and $\frac{1}{2}p_1$ 
 vanishing, namely that coming from 2-torsion.

 \vspace{3mm}
 What we would like to describe is a space, which we will call $\mathcal{F}$, that 
 sits in the fibration
 \(
 \xymatrix{
\mathcal{F}
\ar[rr]  
&&
(BO) \langle 9 \rangle  
\ar[rr]^{\times 8}
&&
(BO) \langle 9 \rangle
}
 \)
 where the second map takes $\frac{p_2}{48}$ to $\frac{p_2}{6}$.
 The question is then:
what is $\mathcal{F}$?
As a warm-up, consider a degree two class, in which case we have 
\(
 \xymatrix{
F
\ar[rr]  
&&
K(\Z, 2) 
\ar[rr]^{\times n}
&&
K(\Z, 2),
}
 \)
 and the answer is $F=K(\Z_n, 1)$. 
 We then seek a lift 
\(
  \xymatrix{
  &&
    K(\Z, 2)
    \ar[dd]^{\times n}
   \\
   \\
   X
   \ar[rr]
   \ar@{..>}[rruu]^{?}
   &&
    K(\Z, 2)
    \ar[rr]^{\alpha_n}
    &&
    K(\Z_n,2),
  }
\)
which exists if and only if the pullback of the map $\alpha_n$ to $X$ is zero,
where  $\alpha_n$ is `reduction mod 2'.
For a String structure on a space $Y$ we have the following diagram
\(
  \xymatrix{
  Y
  \ar[rr]^{x}
\ar[rrdd]
  &&
   \mathcal{F}={\rm String}^{\cal F}
   \ar[rr]
   \ar[dd]
   &&
   K(\Z, 8)
   \ar[dd]^{\times 8}
  &&
   \\
   \\
   &&
   (BO)\langle 8 \rangle
    \ar[rr]_{\frac{1}{6}p_2}
  &&
K(\Z, 8)
\ar[rr]
&& 
K(\Z_8, 8),
  }
  \label{F}
\)
where $x$ is our class $\frac{1}{48}p_2$ which naturally lives 
not in $(BO)\langle 8 \rangle$ but rather in the desired space $\mathcal{F}$.
The above diagram specifies $\mathcal{F}$. The modification of Proposition 
\ref{p26} is then

\begin{proposition}
The class $\frac{1}{48}p_2$ is the obstruction to lifting 
a $String^{\cal F}$ bundle, defined
by diagram (\ref{F}), to a Fivebrane bundle.
\end{proposition}

\vspace{3mm}
We conclude this section with a caveat. 
 The formula that arise from the anomaly
cancelation, e.g. (\ref{dH7}), involve different 
coefficients for each of the second Pontrjagin 
classes of the two bundles. Had we had the same
factor or a just minus sign, then we could have 
simply applied the 
formulae in the discussion at the end of section 
\ref{ineq}. However, we have relative factors of 
$48$, so how can one make sense of such factors?
For instance, in the case when the string condition applies to 
all bundles, is $p_2(E) + \frac{1}{n} p_2(F)$ the
second Pontrjagin class of some combination of 
the bundles $E$ and $F$? It does not make sense
to talk about $E + \frac{1}{n} F$, but one can talk 
about $nE + F$. But then, we would have the issue of
torsion, especially that $n=48$ includes both 
the important 2- and 3-torsion. We can thus say that,
away from such torsion, the Fivebrane 
class of $48TM - E$ is equal to $48$ times the 
Fivebrane class of $TM$ minus the Fivebrane class
of $E$. In the case when $H^8(X, \Z)$ has 
$2$- or 3-torsion we cannot draw such  
a conclusion.

\subsection{Inequivalent Fivebrane Structures: The Fivebrane Class } 
\label{ineq}
In this section we consider the set of inequivalent Fivebrane structures on our manifold or on a pair consisting of a manifold with a gauge bundle,
In order to do so, let us first recall the situation for String structures on manifolds.
The complex version of the fibration 
\(
K(\Z, 3) \to B{\rm String} \to B{\rm Spin} 
\)
is 
\(
K(\Z,3) \to (BU)\langle 6 \rangle \to BSU,
\)
where the space $(BU)\langle 6\rangle$ is the 5-connected cover of $BU$.
A choice
of String structure on a space $X$ is a lift $\hat f$
of the classifying map $f$ in the diagram
\(
  \raisebox{30pt}{
  \xymatrix{
  &&
    (BO)\langle8\rangle
    \ar[dd]
   \\
   \\
   X
   \ar[rr]^{f}
   \ar@{..>}[rruu]^{\hat f}
   &&
    (BO)\langle 4\rangle
    \ar[rr]^{\frac{1}{2}p_1}
    &&
    K(\Z,4).
  }
  }
\label{liftst}
\)
A choice of
$(BU)\langle 6 \rangle$ structure on a space $X$ is a lift in the diagram 
\(
  \raisebox{30pt}{
  \xymatrix{
  &&
    B(U)\langle6\rangle
    \ar[dd]
   \\
   \\
   X
   \ar[rr]^f
   \ar@{..>}[uurr]^{\hat f}
   &&
    BSU
    \ar[rr]^{c_2}
    &&
    K(\Z,4).
  }
  }
  \label{liftcomplex}
\)
The appearance of the factor $\frac{1}{2}$ in (\ref{liftst}) vs. (\ref{liftcomplex}) 
is a reflection of the fact
mentioned earlier that the map from $B{\rm Spin}$ to $BSU$ is given by 
multiplication by two.
The fibration sequences 
\(
\raisebox{15pt}{
\xymatrix@R=4pt{
  K(\Z, 3) \ar[r] & (BO)\langle 8 \rangle \ar[r] & (BO)\langle 4 \rangle 
\\
K(\Z, 3) \ar[r] & (BU)\langle 6 \rangle \ar[r] & BSU 
 }
 }
\)
show that two String or $(BU)\langle 6 \rangle$ structures differ by a map 
to $K(\Z,3)$ for a fixed Spin or $BSU$ structure, respectively.
Therefore,
the set of lifts, i.e. the set of String structures for a fixed Spin structure
in the real case,  or the set of $(BU)\langle 6\rangle$ structures for a fixed 
$BSU$ structure in the complex case, is a torsor
\footnote{Here a torsor is a set with a group action which is free and
transitive. So as a set it is the same as the group. The group acts
but there is no canonical identification of the set with the group.}
 for a quotient
of the third integral cohomology group $H^3(X; \Z)$. The elements of the 
torsor
are the string classes, corresponding to the NS degree three
$H$-field in string theory. 

\vspace{3mm}
Now we are ready to study the set of Fivebrane structures.
The complex version of the fibration 
\(
K(\Z, 7) \to B{\rm Fivebrane} \to B{\rm String}, 
\label{kz71}
\)
obtained by applying the classifying functor on the sequence in Def. \ref{def 2},
is 
\(
K(\Z,7) \to (BU)\langle 10 \rangle \to (BU) \langle 8 \rangle ,
\label{kz72}
\)
where the space $(BU)\langle 8 \rangle$ is the 7-connected cover of $BU$. A choice
of Fivebrane structure on a space $X$ is a lift in the diagram
\(
  \xymatrix{
  &&
    (BO)\langle9\rangle
    \ar[dd]
   \\
   \\
   X
   \ar[rr]^f
   \ar@{..>}[uurr]^{\hat f}
   &&
    (BO)\langle 8\rangle
    \ar[rr]^{\frac{1}{2}p_1}
    &&
    K(\Z,8).
  }
\label{lift5}
\)
A choice of
$(BU)\langle 9 \rangle$ structure on a space $X$ is a lift in the diagram 
\(
  \xymatrix{
  &&
    (BU)\langle 10 \rangle
    \ar[dd]
   \\
   \\
   X
   \ar[rr]^f
   \ar@{..>}[uurr]^{\hat f}
   &&
    (BU)\langle 8  \rangle
    \ar[rr]^{\frac{1}{6}c_4}
    &&
    K(\Z,8).
  }
  \label{liftcomplex}
\)
The fibration sequences (\ref{kz71}) (\ref{kz72})
show that two Fivebrane or $(BU)\langle 9 \rangle$ structures differ by a map 
to $K(\Z,7)$ for a fixed String or $(BU)\langle 7 \rangle$ structure, respectively.

\vspace{3mm}
\noindent {\bf
Remarks.}

{\bf 1.} A degree seven class  $X\to K(Z,7)$ that corresponds to the fivebrane
structure cannot be specified. Consider the diagram
\(
  \raisebox{70pt}{
  \xymatrix{
  &&
  K(\Z, 7)
  \ar[dd]
  \\
  \\
  &&
    (BO)\langle 9\rangle
    \ar[dd]
   \\
   \\
   X
   \ar[rr]^f
   \ar@{..>}[uurr]^{\hat f}
   \ar@{..>}[uuuurr]^{g}
   &&
    (BO)\langle 8\rangle
    &&
  }
  }
\label{3by3}
\)
and note that $f$ and $g$ do not determine ${\hat f}$ since 
$BO\langle 9\rangle \neq BO \langle 8 \rangle \times K(\Z, 7)$. This means
that there is no degree seven class that picks ${\hat f}$. 

\noindent {\bf 2.} For a given string structure, the set of compatible fivebrane structures 
has a transitive action by degree seven cohomology. We have

\begin{proposition}
\label{the set of lifts}
The set of lifts, i.e. the set of Fivebrane structures for a fixed String structure
in the real case,  or the set of $(BU)\langle 9\rangle$ structures for a fixed 
$(BU)\langle 7 \rangle$ structure in the complex case, is a 
torsor for a quotient (to be described below)
of the seventh integral cohomology group $H^7(X; \Z)$. 
\end{proposition}

We call the  elements of the latter the {\it fivebrane classes}, corresponding to the (dual) 
NS degree seven $H$-field in string theory.
We now explain the quotient in the proposition. 
Whether or not we have a 
free action depends on whether the map 
\(
\left[ X,( \Omega (BO)) \langle 8 \rangle \right] \to 
\left[ X, K(Z, 7) \right]
\label{omegakz}
\)
coming from the diagram 
\(
  \raisebox{50pt}{
  \xymatrix{
  &&
 ( \Omega (BO))\langle 8 \rangle
  \ar[d]
\\
  &&
  K(\Z, 7)
  \ar[d]
  \\
  &&
    (BO)\langle9\rangle
    \ar[d]
   \\
   X
   \ar[rr]|f
   \ar@{..>}[urr]|{\hat f}
   \ar@{..>}[uurr]^{g}
   &&
    BO\langle 8\rangle
    &&
  }
  }
\label{omegaBO}
\)
is the zero map.  Note that (\ref{omegakz}) is 
just the map of homotopy classes induced by
the top arrow in (\ref{omegaBO}). However, 
the intention is to consider {\it two} lifitngs ${\hat f}_i$ for 
$i=1,2$ which therefore differ by a map $g$ as in (\ref{3by3}).
Since the  maps of homotopy classes  give long exact sequences 
for any three of the target spaces, then if (\ref{omegakz}) were 
the zero map then the next one would be onto.
Failure to be onto  is measured by the failure of 
$[X, (BO) \langle 8 \rangle] \to \left[ X,K(\Z,7) \right]$ to 
be zero.
The set of lifts is a torsor over the quotient 
\(
H^7(X, \Z) / \left[X, \Omega (BO) \langle 8 \rangle \right].
\)

\noindent {\bf 3.} The fibration (\ref{kz71}) is a map of infinite loop spaces 
so we can think 
of it in terms of a map
\(
K(\Z, 7) \times B{\rm Fivebrane} \to B{\rm Fivebrane},
\)
realizing an action of $K(\Z, 7)$ on the space of inequivalent 
fivebrane structures. 
 
\noindent {\bf 4.} At the level of spectra, which are objects whose homotopy groups 
represent generalized cohomology theories \cite{Adams}, we have
\(
\Sigma^7 H\Z \times ko\langle 9 \rangle \to ko\langle 9 \rangle,
\)
where $ko\langle 9 \rangle$ is the {\it connective} version
of real K-theory $KO$,
i.e. 
the version with no negative degrees,
 of 
the K-theory corresponding to $BO \langle 9 \rangle$.

\noindent {\bf 5.} Addition of vector bundles is encoded in a
product 
\(
BSO \times BSO \to BSO,
\)
which has a lift to 
\(
B{\rm Fivebrane} \times B{\rm Fivebrane} \to B{\rm Fivebrane}.
\label{b555}
\) 
This gives us a way of adding Fivebrane bundles; for example, in 
our setting of heterotic string theory, we have the virtual difference 
of the tangent bundle and the gauge bundle.
In fact lifts such as (\ref{b555}) work for even higher connected covers as well. 

\noindent {\bf 6.} The bundles we have are
the tangent bundle 
and the gauge vector bundle $E$, with the total bundle being the 
K-theoretic virtual difference. Hence it is natural to ask about the 
fivebrane structure of sums and differences of bundles. 
There is a distinction here according to whether the bundles are
spin or not. For the two possible gauge bundles, the $E_8 \times E_8$ 
bundle is Spin, while the ${\rm Spin}(32)/\Z_2$ bundle is not, as its
second Stiefel-Whitney class is non-vanishing. 
The (integral)
Pontrjagin classes for two bundles $E_1$ and $E_2$ satisfy \cite{MilnorStasheff}
\(
p_k(E_1 \oplus E_2) = \sum_{i+j=k} p_i(E_1) \cup p_i(E_2)~~ {\rm mod~ elements ~of~ order~} 2.
\)
Thus, applying to our case,
\(
p_2(TM \oplus E) = p_1(TM) \cup  p_1(E) + p_2(TM) + p_2(E)~~{\rm mod~} 2-{\rm torsion}.
\)
But we are assuming $p_1(TM)$ to be zero -- recall that in general, we have 
to have a String structure on a bundle
to talk about Fivebrane structure as we have to have a Spin structure
on a bundle to talk about 
String structure. Hence we get:
\(
p_2(TM \oplus E)= p_2(TM) + p_2(E)~~ {\rm mod~}2-{\rm torsion},
\)
so
\(
2[ p_2(TM \oplus E) - p_2(TM) - p_2(E)]=0.
\)
We have worked out the case of a direct sum, but the case of virtual difference
should be analogous. Indeed, since $V \oplus (-V) $ is equivalent to a trivial bundle, 
 $p_1(E_1-E_2)=p_1(E_1) - p_1(E_2)$ and 
$p_2(E_1-E_2)=p_2(E_1) - p_2(E_2)$, and hence we get
\(
2[ p_2(TM - E) - p_2(TM) + p_2(E)]=0.
\)
We can actually give a better description by using special
characteristic classes 
more adapted for covers of the tangent bundle and, 
in addition, use integral coefficients. 

\vspace{3mm}
Note that,
making use of the fact that the bundles are not just real but
also spin, 
\( 
H^*(B{\rm Spin}; \Z)=\Z[Q_1,Q_2,\cdots]\oplus \gamma,
\label{qs}
\) 
with $\gamma$ a 2-torsion factor, i.e. $2\gamma=0$ \cite{Thomas},
concentrated in degrees \emph{not} congruent to 0 mod 4.
The two degrees relevant to our discussion are
 \bea H^4(B{\rm Spin};
\Z)&\cong&\Z ~~~~~~~~~~~{\rm with~~generator}~~~ Q_1
\nonumber\\
H^8(B{\rm Spin}; \Z)&\cong&\Z\oplus \Z ~~~~~{\rm with~~generators}~~~ Q_1^2,
Q_2,
\eea
where $Q_1$ and $Q_2$ are determined by their relation to the Pontrjagin classes
\begin{eqnarray}
p_1&=&2Q_1 \nonumber\\
p_2&=&Q_1^2 + 2Q_2.
\label{Qs}
\end{eqnarray}
Then the spin generators are given in terms of the Pontrjagin classes by
$Q_1=\frac{1}{2}p_1$ and $Q_2=\frac{1}{2}p_2 - \frac{1}{2}\left(\frac{1}{2}p_1\right)^2$. 
Note that
$Q_2=\frac{1}{2}p_2$ holds when $Q_1$  vanishes.
Their importance in the study of anomalies was emphasized in \cite{SSpin}. 
For two Spin bundles $E$ and $F$, such as the spinor bundle and the 
$E_8 \times E_8$ bundle
\bea
Q_2(E \oplus F) &=& \sum_{i+j= 2}Q_i(E) \cup Q_j(F)
\nonumber\\
&=& Q_2(E) + Q_2(F),
\label{add}
\eea
where again the condition $Q_1(TM)=0$ is used.

\subsection{Overview}

\begin{itemize}

  \item
    We considered the notion of a ``Fivebrane structure'' on a manifold,
    which generalizes that of a String structure: as a String structure is
    defined to be a lift of the tangent bundle of a Spin-manifold to 
    the 3-connected cover of Spin, a Fivebrane structure is the further lift
    to a 7-connected cover. We showed that Fivebrane structures exist when 
    a fractional second Pontrjagin class vanishes. This holds for the 
    tangent bundle as well as gauge bundles. 
    When considering direct sums or K-theoretic 
    virtual bundles, as in the case of heterotic string theory where we have both
    $TX$ and $E$, the obstructions are given by the sums
    and differences of the individual fivebrane structures,
    respectively.

\item  In analogy to how String structures appear in terms of the vanishing of a worldsheet 
anomaly for the superstring, Fivebrane structures are related to
the vanishing of an anomaly in the 6-dimensional worldvolume theory of 
the fivebranes. In M-theory and type IIA string theory, this is the 
anomaly of the M-theory fivebrane \cite{Dixoncoupling}    
\cite{Dixonchern}
\cite{Percacci}
\cite{Bergshoeff1}
\cite{Cederwall}, while in heterotic string theory this is the worldvolume 
anomaly of the heterotic fivebrane \cite{LT}. Note that the electric-magnetic duality
between Strings and 5-branes in ten dimensions relates String structures and
Fivebrane structures:

  \item  We notice that string theory suggests that String structures and Fivebrane 
    structures 
    are related
    by a duality which generalizes Hodge duality.
    Apart from the 
    well-known fact (see, for instance, \cite{DKL} for a survey) that 
    NS 5-branes are magnetic duals to strings, we notice
    that there is a known formula for the Hodge dual of the 3-form curvature
    which appears in the Green-Schwarz mechanism which relates to Fivebrane
    structures as the former relates to String structures.
We leave the detailed study of this duality for a separate treatment. 

\end{itemize}

\appendix

\section{Recollection of characteristic classes}

We recall elements of the theory of characteristic classes.

\subsection{Universal characteristic classes}
For any topological group $G$, there is a \emph{classifying space} $BG$ and a
universal principal $G$-bundle 
\(
EG\to BG
\)
such that equivalence classes of (numerable) $G$-bundles $E\to B$ are in 1-1 correspondence with homotopy classes of maps $B\to BG.$

\vspace{3mm}
For $G$ a Lie group, elements of the cohomology $H^\bullet(B G)$ (with any coefficients)
are called the \emph{universal characteristic classes for $G$-bundles}.
Given a $G$-bundle $E\to B$, its characteristic classes are the pull-backs
of the universal characteristic classes via the classifying map $B\to BG.$
It is a classical theorem that, for $G$ compact  connected,
the image in real cohomology, $H^\bullet(BG, \mathbb{R})$,
 of 
this cohomology ring is finitely freely generated in even degree and 
is isomorphic to the the ring of invariant polynomials
on $\gg = \mathrm{Lie}(G)$. Moreover, the real cohomology ring of $G$ itself
$H^\bullet(G,\mathbb{R}) \simeq H^\bullet(\gg,\mathbb{R})$
is isomorphic to the Lie algebra cohomology ring of $\gg$, which is
generated in odd degree.
The relation between $H^*(G, \mathbb{R})$ and  $H^*(BG ,\mathbb{R})$
is as follows: 

\vspace{2mm}
$\bullet$ $H^*(G)$ is isomorphic to an exterior algebra $\Lambda (x_1,\cdots, x_n)$,
where the $x_i$ are all of odd degree.

$\bullet$ $H^*(BG)$ is isomorphic to a polynomial algebra 
$P[y_1,\cdots, y_n]$ where 
\(
{\rm degree}\  y_i = {\rm degree}\  x_i +1.
\)
The isomorphism of the respective vector spaces of generators can in fact be expressed invariantly:
Let $\mathcal{P}H^*(G)$ denote the subspace of \emph{primitive} elements, i.e. those $x\in H^*(G)$ 
such that 
\(
m^*(x) = 1\otimes x + x\otimes 1 
\)
where $m: G\times G\to G$ is the group multiplication.
Let $\mathcal{Q}H^*(BG)$ denote the \emph{quotient space of indecomposables}, i.e. $H^+(BG)/H^+(BG)\cdot H^+(BG)$ where $H^+$ denotes the subalgebra of elements of positive degree. The isomorphism we have described in terms of generators is in fact 
\(
\tau:\mathcal{P}H^*(G)\to\mathcal{Q}H^{*+1}(BG).  
\)
The isomorphism says that the two vector spaces are
\emph{transgressively related}.  Unfortunately the name \emph{transgression} is sometimes applied
in algebraic topology
to $\tau$  (cf. H. Cartan, Borel, Serre) and sometimes in differential geometry to its inverse 
(cf. Chern and disciples). In algebraic topology, the map $\sigma: H^*(BG) \to H^{*-1}(G)$
is always defined whereas $\tau$ is defined only on primitive elements.

\vspace{3mm}
 An element $c$ in $H^n(BG ,\mathbb{Z})$ can be identified with a 
(homotopy class of a) map
$f$ from $B G$ to the Eilenberg-MacLane space $K(\mathbb{Z},n)$ 
(and similarly for $\mathbb{R}$),
then the 
characteristic class $f^* c$ of a $G$-bundle $E$ corresponding to the 
characteristic class $c$ of $G$ is given simply by the composite map
\(
  \xymatrix{
    B \ar[rr]^f
    &&
    B G
    \ar[rr]^c
    &&
    K(\mathbb{Z}, n)
  }
  \,.
\)

\subsection{The cohomology of $BSO$ and $BU$}
We are particularly concerned in this paper with $H^*(BSO)$ and $H^*(BU).$ 
See \cite{MilnorStasheff}.
In characteristic 0, 
\(
H^*(BSO)\simeq P[y_4,y_8,\cdots],
\)
 where now the indexing is by the
degree of the (Pontrjagin) class.
In any characteristic, 
\(
H^*(BU)\simeq P[c_1,c_2,\cdots],
\)
 where we have used the traditonal notation for the Chern classes $c_i$.
The basic definitions and properties of 
Chern and Pontrjagin classes
are recalled next.
 
 \subsection{Chern-Weil homomorphism}


Given a $G$-bundle $P \to X$ with connection
$A$, the
Chern-Weil homomorphism
\(
  \mathrm{inv}(\gg) \to H^\bullet(X,\mathbb{R})
\)
sends each degree $n$ invariant polynomial $k \in \mathrm{inv}(\gg)$ 
on the Lie algebra $\gg = \mathrm{Lie}(G)$ to the differential form
\(
  k \mapsto k(F_A) := k(F_A \wedge \cdots \wedge F_A)
\)
obtained by wedging $n$ copies of the curvature 2-form and 
evaluating its Lie algebra value in $k$. This $k(F_A)$ is a closed
form. The corresponding class $[k(F_A)]$ in deRham cohomology is the
characteristic class of $P$ corresponding to $k$. This class
is independent of the choice of connection on $P$.

\begin{table}[h]
\begin{center}
$$
\xymatrix{
   \mathrm{inv}(\gg)
   \ar[rr]
   \ar@/_2pc/[rrrr]|{\mbox{Chern-Weil homomorphism}}
   &&
   \Omega^\bullet(X)  
   \ar[rr]
   && 
   H^\bullet(X,\mathbb{R}) 
   \\
   k
   \ar@{|->}[rr] 
   && 
   k(F_A) 
   \ar@{|->}[rr] 
   &&  
   [k(F_A)] 
   \\
   \mbox{
   \begin{tabular}{c}
     invariant
     \\
     polynomial
   \end{tabular} 
   }
   &
   &
   \mbox{
   \begin{tabular}{c}
     characteristic
     \\
     form
   \end{tabular}
   }
   &&
   \mbox{
   \begin{tabular}{c}
     characteristic
     \\
     class
   \end{tabular}
   }
}
$$
\end{center}
\caption{
  {\bf The Chern-Weil homomorphism} sends, for each $G$-bundle $P \to X$,
  any degree $n$ invariant polynomial on $\gg = \mathrm{Lie}(G)$ to the 
  deRham class of the differential form $k(F_A) = k(F_A \wedge \cdots \wedge F_A)$
  obtained by inserting the curvature 2-form of any connection on $P$ into $k$.
}
\end{table}

\subsection{Polynomials and classes for matrix Lie algebras}

Given a matrix Lie algebra $\gg\subset \mathfrak{gl}(n)$, the trace and the
determinant operation on matrices provide families of $\gg$-invariant polynomials.

\begin{itemize}

\item

The assignment
\(
  X \mapsto \mathrm{ch}(X) := \mathrm{Tr}(\exp(X))
\)
defines the invariant polynomials $\mathrm{ch}_k(X)$
as
\(
  \mathrm{ch}(X) := \sum\limits_{k=0}^\infty t^k \mathrm{ch}_k(X,\cdots,X)
  \,.
\)

The expression $\mathrm{ch}(X)$ is the total Chern-character. 

\item

Let $\gg \subset \mathfrak{gl}(c,\mathbb{C})$ be a Lie algebra of
\emph{complex} matrices. Then the assignment
\(
  X \mapsto c(X) := \mathrm{det}(t + i X)
\)
defines the invariant polynomials $c_k$ as
\(
  c(X) = \sum\limits_{k=0}^{n} t^{n-k} c_k(X,\cdots,X)
  \,.
\)
These are the Chern Polynomials. The corresponding class $[c_k(F_A)]$
is the $k$th Chern class.

\item

Let $\gg \subset \mathfrak{gl}(n,\mathbb{R})$  be a Lie algebra of
\emph{real} matrices. Then the assignment
\(
  X \mapsto c(X) := \mathrm{det}(t - X)
\)
defines the invariant polynomials $p_{k/2}$ as
\(
  c(X) = \sum\limits_{k=0}^{n} t^{n-k} p_{k/2}(X,\cdots,X)
  \,.
\)
These are the Pontrjagin polynomials. The corresponding classes  $[p_{k/2}(F_A)]$
are the Pontrjagin classes.

\item

  The restrictions of the $c_k$ from $\mathfrak{gl}(n,\mathbb{C})$
  to $\mathfrak{gl}(n,\mathbb{R})$ satisfy
  \(
    i^k c_k(X,\cdots,X) = p_{k/2}(X,\cdots, X)
    \,.
  \)

\item

  When $\gg \subset \mathfrak{sl}(n)$ all elements are traceless and
  various cancellations occur. In particular for 
  $\gg = \mathfrak{so}(n)$ the Pontrjagin classes have
  relatively simple relation to the Chern characters.

\end{itemize}

\begin{table}[h]
 \begin{center}
 \begin{tabular}{c|cccc}
   {\bf Lie algebra}
   &
   $\gg \subset \mathfrak{gl}(n,\mathbb{C})$
   &
   $\gg \subset \mathfrak{gl}(n,\mathbb{C})$   
   &
   $\gg \subset \mathfrak{gl}(n,\mathbb{R})$   
   \\
   \hline
   \\
   \begin{tabular}{l}
     {\bf invariant }
     \\
     {\bf polynomial}
   \end{tabular}
   &
   $\mathrm{ch}_k$
   &
   $\mathrm{c}_k$
   &
   $\mathrm{p}_{k/2}$
   \\
   \hline
   \\
   {\bf definition}
   &
   \begin{tabular}{l}
     $\mathrm{ch}(X)$\\ 
     $= \mathrm{tr}(\exp(i t X))$
     \\
     $=\sum\limits_k t^k \mathrm{ch}_k(X,\cdots,X)$
   \end{tabular}
   &
   \begin{tabular}{l}
     $\mathrm{c}(X)$\\ 
     $= \mathrm{det}(t + i X)$
     \\
     $=\sum\limits_k t^{n-k} \mathrm{c}_k(X,\cdots,X)$
   \end{tabular}
   &
   \begin{tabular}{l}
     \\ 
     $\mathrm{det}(t -  X)$
     \\
     $=\sum\limits_k t^{n-k} p_{k/2}(X,\cdots,X)$
   \end{tabular}
   \\
   \hline
   \\
   \begin{tabular}{l}
     {\bf characteristic}
     \\
     {\bf class}
   \end{tabular}
   &
   \begin{tabular}{c}
     Chern character
   \end{tabular}
   &
   \begin{tabular}{c}
     Chern class
   \end{tabular}
   &
   \begin{tabular}{c}
     Pontrjagin class
   \end{tabular}
   \\
   \hline
 \end{tabular}
 \end{center}
\vspace{2mm} \caption{{\bf 
   Characteristic classes for matrix Lie algebras}
   obtained from the trace and the determinant.
 }
\end{table}

\section{The Leray-Serre spectral sequence}
Suppose that 
$H^*\left((BU)\langle 2k \rangle; \Q \right)= P \left[ c_k, c_{k+1}, \cdots \right]$.
Corresponding to the fibration
\(
  \xymatrix{
  K(\Z, 2k-1)
  \ar[rr]
  &&
    (BU)\langle 2k+2\rangle
    \ar[dd]
    &&
   \\
   \\
   &&
   (BU)\langle 2k \rangle
   \ar[rr]
    &&
    K(\Z,2k)
  }
  \label{SerreSS}
\)
we have 
\(
H^* \left( (BU)\langle 2k \rangle; H^* \left( K(\Z, 2k-1);\Q \right) \right)
\Longrightarrow H^*\left( (BU)\langle 2k+2 \rangle \right).
\)
Here the shift in 2 in the degree of the Eilenberg-MacLane space 
comes from Bott periodicity. The cohomology 
$H^* \left( K(\Z, 2k-1);\Q \right)$ is an exterior algebra $\Lambda(e_{2k-1})$
on a generator of degree $2k-1$, $e_{2k-1}$, which satisfies $e_{2k-1}^2=0$. 
The differential in the sequence takes
this generator to the $k$-th Chern class $c_k$. Modding out by this class 
gives (\ref{rationalBU}). In fact, the spectral sequence collapses from here and there is no extension problem for the $E_\infty$ term.
  
  \vspace{3mm}
   Since $G$ is 
simply connected, the $E_2$ term of this spectral sequence is in 
bidegree $(p,q)$ 
\( 
E_2^{p,q} = H^p(G)\otimes H^q(K(\Z,2) 
\) 
The cohomology of $K(\Z,2)$ is a polynomial algebra on one 
generator $y$ in degree $2$: $H(K(\Z,2))= \R[y]$, and the differential 
$d_2\colon E_2^{p,q} \to E_2^{p+2,q-1}$ raises total degree by 
$1$.  
\( 
E_2^{p,q} = \begin{cases} 
H^p(G)\otimes \R y^{k}\ \text{if}\ q = 2k\ \text{is even} \\ 
0 \ \text{if}\ q\ \text{is odd.} 
\end{cases} 
\)
Moreover, since we are dealing with real coefficients, the cohomology 
of $G$ is an exterior algebra $H(G) = \Lambda (x_3,x_5,\ldots )$.  So 
we are computing the cohomology of a complex which looks like  
\( 
\Lambda V\otimes \R[y],
\) 
where $E$ is a graded vector space concentrated in odd degrees and 
$x$ is of degree $2$.  Clearly the differential $d_2$ is zero so that 
$E_3 = E_2$.  The new differential $d_3\colon E_3^{p,q}\to E_3^{p+3,q-2}$
 sends $y$ to $x_3 = H$, the 
generator of $H^3(G)$.  Since $d_3$ is 
a derivation with 
respect to the algebra structure of $E_3^{p,q}$. It follows that 
\( 
d_3(z\otimes y^k) = k Hz\otimes y^{k-1} 
\) 
where $z\in H(G)$. 

  \bigskip\bigskip
\noindent
{\bf \large Acknowledgements}\\
\vspace{1mm}
\noindent

H.S. thanks Matthew Ando and Nitu Kitchloo for very useful 
discussions and acknowledges the hospitality of the mathematics 
department at the University of Illinois at Urbana-Champaign. 
We thank Dan Freed and Jacques Distler for helpful discussion.
  

\vspace{1cm}

{\small

\begin{tabular}{ccc}

\begin{tabular}{l}
Hisham Sati\\
Department of Mathematics\\Yale University\\10 Hillhouse Avenue\\
New Haven, \\
CT 06511
\end{tabular}
&
\begin{tabular}{l}
Urs Schreiber\\
Fachbereich Mathematik\\Schwerpunkt Algebra und Zahlentheorie\\Universit\"at
Hamburg\\Bundesstra\ss e 55\\D--20146 Hamburg 
\end{tabular}
&
\begin{tabular}{l}
Jim Stasheff\\
Department of Mathematics\\University of Pennsylvania\\
David Rittenhouse Lab. \\
209 South 33rd Street \\
Philadelphia, PA 19104-6395
\end{tabular}

\end{tabular}
}

 \end{document}